\documentstyle[12pt]{article}
\input amssym.def
\parskip=\medskipamount
\newtheorem{definition}{Definition}[section]

\newtheorem{theorem}[definition]{Theorem}
\newtheorem{proposition}[definition]{Proposition}
\newtheorem{corollary}[definition]{Corollary}
\newtheorem{remark}[definition]{Remark}

\newtheorem{examples}[definition]{Examples}

 1   
\font\ddpp=msbm10  scaled \magstep 1  

\def\QED{\hskip0.1em\hfill\null\ \null\nobreak\hfill
\kern3pt\lower1.8pt\vbox{\hrule\hbox
{\vrule\kern1pt\vbox{\kern1.7pt \hbox{$\scriptstyle
QED$}\kern0.2pt}\kern1pt\vrule}\hrule}}
\def\R{\hbox{\ddpp R}}               
\def\C{\hbox{\ddpp C}}    

\def\lcf{\lbrack\! \lbrack}
\def\rcf{\rbrack\! \rbrack}

\newcommand\prueba {\mbox{{\em Proof: }}}

\begin{document}
\baselineskip=.55cm
\title{{\bf GENERALIZED LIE BIALGEBROIDS AND JACOBI STRUCTURES}}
\author{David Iglesias, Juan C. Marrero
\\ {\small\it Departamento de Matem\'atica
Fundamental, Facultad de Matem\'aticas,}\\[-8pt] {\small\it
Universidad de la Laguna, La Laguna,} \\[-8pt] {\small\it
Tenerife, Canary Islands, SPAIN,}\\[-8pt] {\small\it E-mail:
diglesia@ull.es, jcmarrer@ull.es} }
\date{}

\maketitle
\baselineskip=.3cm
\begin{abstract}
{\small
The notion of a generalized Lie bialgebroid (a generalization of the
notion of a Lie bialgebroid) is introduced in such a way that a Jacobi
manifold has associated a canonical generalized Lie bialgebroid. As a
kind of converse, we prove that a Jacobi structure can be defined on
the base space of a generalized Lie bialgebroid. We also show that it
is possible to construct a Lie bialgebroid from a generalized Lie
bialgebroid and, as a consequence, we deduce a duality theorem.
Finally, some special classes of generalized Lie bialgebroids are
considered: triangular generalized Lie bialgebroids and generalized
Lie bialgebras.
}
\end{abstract}
\begin{quote}

\end{quote}

\section{Introduction}
\baselineskip=.55cm
Roughly speaking, a Lie algebroid over a manifold $M$ is a vector
bundle $A$ over $M$ such that its space of sections $\Gamma (A)$
admits a Lie algebra structure $\lcf \, ,\, \rcf$ and, moreover,
there exists a bundle map $\rho$ from $A$ to $TM$ which provides a Lie
algebra homomorphism from $(\Gamma (A),\lcf \, ,\, \rcf )$ into the
Lie algebra of vector fields $\frak X (M)$ (see \cite{Mk,Pr}). Lie
algebroids are a natural generalization of tangent bundles and real Lie
algebras of finite dimension. But, there are many
other interesting examples, for instance, the cotangent bundle
$T^\ast M$ of any Poisson manifold  $M$ possesses a
natural Lie algebroid structure (\cite{BV,CDW,F,V}). In fact, there
is a one-to-one correspondence between Lie algebroid structures on a
vector bundle $A$ and linear Poisson structures on the dual bundle
$A^\ast$ (see \cite{CDW,Co}). An important class of Lie algebroids
are the so-called Lie bialgebroids. This is a Lie algebroid $A$ such
that the dual vector bundle $A^\ast$ also carries a Lie algebroid
structure which is compatible in a certain way with that on $A$ (see
\cite{K-S,MX}). If $M$ is a Poisson manifold, then the pair $(TM,
T^\ast M)$ is a Lie bialgebroid. As a kind of converse, it was proved
in \cite{MX} that the base space of a Lie bialgebroid is a Poisson
manifold. Apart from the pair $(TM,T^\ast M)$ ($M$ being a Poisson
manifold), other interesting examples of Lie bialgebroids are Lie
bialgebras \cite{D}. A Lie bialgebra is a Lie bialgebroid such that
the base space is a single point and there is a one-to-one
correspondence between Lie bialgebras and connected simply connected
Poisson Lie groups (see \cite{D,LW,V}). We remark that a connected
simply connected abelian Poisson Lie group is isomorphic to the dual
space of a real Lie algebra endowed with the usual linear Poisson
structure (the Lie-Poisson structure). Moreover, Poisson Lie groups
are closely related with quantum groups (see \cite{D2}).

As it is well-known, a Jacobi structure on a manifold $M$ is a
2-vector $\Lambda$ and a vector field $E$ on $M$ such that $[\Lambda
,\Lambda ]=2E\wedge \Lambda$ and $[E,\Lambda ]=0$, where $[\, ,\, ]$
is the Schouten-Nijenhuis bracket \cite{Li2}. If $(M,\Lambda ,E)$ is
a Jacobi manifold one can define a bracket of functions, the Jacobi
bracket, in such a way that the space $C^\infty (M,\R )$ endowed with
the Jacobi bracket is a local Lie algebra in the sense of Kirillov
\cite{K}. Conversely, a local Lie algebra structure on $C^\infty
(M,\R )$ induces a Jacobi structure on $M$ \cite{GL,K}. Jacobi
manifolds are natural generalizations of Poisson manifolds. However,
very interesting manifolds like contact and locally conformal
symplectic (l.c.s.) manifolds are also Jacobi and they are not
Poisson. In fact, a Jacobi manifold admits a generalized foliation
whose leaves are contact or l.c.s. manifolds (see \cite{DLM,GL,K}).
If $M$ is an arbitrary manifold, the vector bundle $TM\times \R \to
M$ possesses a natural Lie algebroid structure. Moreover, if $M$ is a
Jacobi manifold then the 1-jet bundle $T^\ast M\times \R \to M$
admits a Lie algebroid structure \cite{KS} (for a Jacobi manifold the
vector bundle $T^\ast M$ is not, in general, a Lie algebroid).
However, the pair $(TM\times \R ,T^\ast M\times \R )$ is not, in
general, a Lie bialgebroid (see \cite{V2}).

On the other hand, in \cite{IM}, we studied Jacobi structures on the dual
bundle $A^\ast$ to a vector bundle $A$ such that the Jacobi bracket
of linear functions is again linear and the Jacobi bracket of a
linear function and the constant function 1 is a basic function. We
proved that a Lie algebroid structure on $A$ and a 1-cocycle $\phi _0\in \Gamma
(A^\ast )$ induce a Jacobi structure on $A^\ast$ which satisfies the above
conditions. Moreover, we showed that this correspondence is a
bijection. We also consider two interesting examples: i) for an arbitrary
manifold $M$, the Lie algebroid $A=TM\times \R$ and the 1-cocycle
$\phi _0=(0,1)\in \Omega ^1(M)\times
C^\infty (M,\R )\cong\Gamma (A^\ast )$, we prove that the resultant
linear Jacobi structure on $T^\ast M\times \R$ is just the canonical
contact structure and ii) for a Jacobi manifold $M$, the Lie
algebroid $A^\ast = T^\ast
M\times \R$ and the 1-cocycle $X_0=(-E,0)\in \frak X (M)\times
C^\infty (M,\R )\cong \Gamma (A)$, we deduce that the corresponding linear
Jacobi structure $(\Lambda_{(TM\times \R ,X_0)},E_{(TM\times \R
,X_0)})$ on $TM\times \R$ is given by
$$\Lambda_{(TM\times \R ,X_0)}=\Lambda ^c+\frac{\partial}{\partial
t} \wedge E^c-t\Big( \Lambda ^v+\frac{\partial}{\partial t}\wedge
E^v\Big) ,\qquad E_{(TM\times \R ,X_0)}=E^v,$$
where $\Lambda ^c$ (resp. $\Lambda^v$) is the complete (resp. vertical)
lift to $TM$ of $\Lambda$ and $E^c$ (resp. $E^v$) is the complete
(resp. vertical) lift to $TM$ of $E$. This Jacobi structure was
introduced in \cite{ILMM} and it is the Jacobi counterpart to the
tangent Poisson structure first used in \cite{S} (see also \cite{Co,GU}).

Therefore, for a Jacobi manifold $M$, it seems reasonable to consider
the pair $((A=TM\times \R ,\phi _0=(0,1))$,$(A^\ast=T^\ast M
\times \R ,X_0=(-E,0)))$ instead of the pair $(TM\times \R ,T^\ast M
\times \R )$. In fact, we prove, in this paper, that the Lie
algebroids $TM\times \R$ and $T^\ast M\times \R$ and the 1-cocycles
$\phi _0$ and $X_0$ satisfy some compatibility conditions. These
results suggest us to introduce, in a natural way, the definition of
a generalized Lie bialgebroid. The aim of this paper is to discuss
some relations between generalized Lie bialgebroids and Jacobi structures.

The paper is organized as follows. In Section 2, we recall several
definitions and results about Jacobi manifolds and Lie algebroids
which will be used in the sequel. In Section 3, we show some facts
about the differential calculus on Lie algebroids in the presence of
a 1-cocycle. If $(A,\lcf \, ,\, \rcf ,\rho )$ is a Lie algebroid over
$M$ and, in addition, we have a 1-cocycle $\phi _0\in \Gamma (A
^\ast)$ then the usual representation of the Lie algebra $\Gamma (A)$
on the space $C^\infty (M,\R )$ can be modified and a new
representation is obtained. The resultant cohomology operator $d_{\phi
_0}$ is called the $\phi _0$-differential of $A$ and its expression,
in terms of the differential $d$ of $A$, is $d_{\phi
_0}\omega =d\omega +\phi _0\wedge \omega$, for $\omega \in \Gamma
(\wedge ^k A^\ast )$. The $\phi _0$-differential of $A$ allows us to
define, in a natural way, the $\phi _0$-Lie derivative by a section
$X\in \Gamma (A)$, $({\cal L}_{\phi _0})_X$, as the commutator of
$d_{\phi _0}$ and the contraction by $X$, that is, $({\cal L}_{\phi
_0})_X  =d_{\phi _0}\circ i_X+i_X\circ d_{\phi _0}$. On the other
hand, imitating the definition of the Schouten bracket of two
multilinear first-order differential operators on the space of
$C^\infty$ real-valued functions on a manifold $N$ (see \cite{BV}),
we introduce the $\phi _0$-Schouten bracket of a $k$-section $P$ and a
$k'$-section $P'$ as the $k+k'-1$-section given by
$$
\lcf P, P' \rcf_{\phi _0}=\lcf P, P' \rcf + (-1)^{k+1}(k-1)P\wedge
(i_{\phi _0}P')-(k'-1)(i_{\phi_0}P)\wedge P',
$$
where $\lcf \, ,\, \rcf$ is the usual Schouten bracket of $A$ (for
the properties of the $\phi _0$-Schouten bracket, see Theorem
\ref{teoremaschou}).  When $(M,\Lambda ,E)$ is a Jacobi manifold and
we consider the Lie algebroids $TM\times \R$ and $T^\ast M\times \R$
and the 1-cocycles $\phi _0=(0,1)$ and $X_0=(-E,0)$, we prove that
the above operators satisfy certain compatibility conditions (see
Proposition \ref{ejemplos}). These results suggest us to introduce, in
Section 4, the definition of a generalized Lie bialgebroid as
follows. Suppose that $(A,\lcf \, ,\, \rcf ,\rho)$ is a Lie algebroid and
that $\phi _0\in \Gamma (A^\ast )$ is a 1-cocycle. Assume also that
the dual bundle $A^\ast$ is a Lie algebroid and that $X_0\in \Gamma
(A)$ is a 1-cocycle. Then, the pair $((A,\phi _0),(A^\ast ,X_0))$ is
said to be a generalized Lie bialgebroid over $M$ if for all $X ,Y
\in \Gamma (A)$ and $P\in \Gamma (\wedge ^kA)$
$$
 d_\ast{}_{X_0}\lcf X ,Y \rcf = \lcf X ,d_\ast{}_{X_0}Y \rcf _{\phi
_0}-\lcf Y ,d_\ast{}_{X_0}X \rcf _{\phi _0},
$$
$$
 ({\cal L}_\ast {}_{X_0})_{\phi _0}P+\lcf X_0,P\rcf _{\phi _0}= 0,
$$
where $d_\ast{}_{X_0}$ and ${\cal L}_\ast{}_{X_0}$ are the
$X_0$-differential and the $X_0$-Lie derivative, respectively, of
$A^\ast$. If the 1-cocycles $\phi _0$ and $X_0$ are null, we recover
the notion of a Lie bialgebroid. Moreover, if $(M,\Lambda ,E)$ is a
Jacobi manifold, the pair $((TM\times \R ,(0,1))$,$(T^\ast M \times \R
,(-E,0)))$ is a generalized Lie bialgebroid. In fact, extending a
result of \cite{MX} for Lie bialgebroids, we prove that the base
space of a generalized Lie bialgebroid is a Jacobi manifold (Theorem
\ref{unteorema}). It is well-known that the product of a Jacobi
manifold with $\R$,
endowed with the Poissonization of the Jacobi structure, is a Poisson
manifold (see \cite {Li2} and Section \ref{introduccion}). We show a
similar result for generalized Lie bialgebroids. Namely, we prove
that if $((A,\phi _0),(A^\ast ,X_0))$ is a generalized Lie
bialgebroid over $M$ then it is possible to define a Lie bialgebroid
structure on the dual pair of vector bundles $(A\times \R ,A^\ast
\times \R)$ over $M\times \R$, in such a way that the induced
Poisson structure on $M\times \R$ is just the Poissonization of the
Jacobi structure on $M$ (Theorem \ref{bialgebrizacion}). Using this
result, we deduce that the generalized Lie bialgebroids satisfy a
duality theorem, that is, if $((A,\phi _0),(A^\ast ,X_0))$ is a
generalized Lie bialgebroid, so is $((A^\ast ,X _0),(A,\phi _0))$.

In Section 5, we prove that it is possible to obtain a generalized Lie
bialgebroid from a Lie algebroid $(A,\lcf \, ,\, \rcf ,\rho )$, a
1-cocycle $\phi _0$ on it and a bisection $P\in \Gamma (\wedge ^2A)$
satisfying $\lcf P,P\rcf _{\phi _0}=0$ (Theorem \ref{triangulares}).
This type of generalized Lie bialgebroids are called triangular.
Examples of triangular generalized Lie bialgebroids are the
triangular Lie bialgebroids in the sense of \cite{MX} and the
generalized Lie bialgebroid associated with a Jacobi structure.
Finally, in Section 6, we study generalized Lie bialgebras, i.e.,
generalized Lie bialgebroids over a single point. Using the results of
Section 5, we deduce that generalized Lie bialgebras can be obtained
from algebraic Jacobi structures on a Lie algebra. This fact allows
us to give examples of Lie groups whose Lie algebras are generalized
Lie bialgebras. The study of this type of Lie groups is the subject
of a forthcoming paper \cite{IM2}.
\section{Jacobi manifolds and Lie algebroids}
\setcounter{equation}{0}
Let $M$ be a differentiable manifold of dimension $n$. We will denote
by $C^\infty (M,\R)$ the algebra of $C^\infty$ real-valued functions
on $M$, by $\Omega ^k(M)$ the space of $k$-forms, by ${\cal V}^k (M)$
the space of $k$-vectors, with $k\geq 2$, by $\frak X (M)$ the Lie
algebra of vector fields, by $\delta$ the usual differential on
$\Omega ^\ast (M)=\oplus _k \Omega ^k (M)$ and by $[\, ,\, ]$ the
Schouten-Nijenhuis bracket (\cite{BV,V}).
\subsection{Jacobi manifolds}\label{introduccion}
A {\em Jacobi structure} on $M$ is a pair $(\Lambda ,E)$, where $\Lambda$
is a 2-vector and $E$ is a vector field on $M$ satisfying the following
properties:
\begin{equation}\label{ecuaciones}
              [\Lambda ,\Lambda ]=2E\wedge \Lambda ,\hspace{1cm}
               [E,\Lambda ]=0.
\end{equation}
The manifold $M$ endowed with a Jacobi structure is called a {\em Jacobi
manifold}. A bracket of functions (the {\em Jacobi bracket}) is defined by
              $$\{ f,g\} =\Lambda (\delta f,\delta g)+fE(g)-gE(f),$$
for all $f,g\in C^\infty (M,\R )$. In fact, the space $C^\infty
(M,\R)$ endowed with the Jacobi bracket is a
{\em local Lie algebra} in the sense of Kirillov (see \cite{K}), that
is, the mapping $\{ \, ,\, \}:C^\infty (M,\R)\times C^\infty
(M,\R)\to C^\infty (M,\R)$ is $\R$-bilinear, skew-symmetric, satisfies the
Jacobi identity and
$$support \{ f ,g\}\subseteq support \, f\cap support \, g ,\mbox{ for
}f,g \in C^\infty (M,\R) ,$$
or, equivalently, $\{ \, ,\, \}$ is $\R$-bilinear, skew-symmetric,
satisfies the Jacobi identity and is a first-order differential operator
on each of its arguments, with respect to the ordinary multiplication
of functions, i.e.,
$$\{ f_1f_2,g\} =f_1\{ f_2,g\}+f_2\{ f_1,g\} - f_1f_2\{ 1,g\} ,\mbox{
for }f_1,f_2,g \in C^\infty (M,\R).$$
Conversely, a structure of local Lie algebra on $C^\infty (M,\R)$ defines a
Jacobi structure on $M$ (see \cite{GL,K}). If the vector field $E$ identically
vanishes then $(M,\Lambda )$ is a {\em Poisson manifold}. Jacobi and Poisson
manifolds were introduced by Lichnerowicz (\cite{Li1,Li2}) (see also
\cite{BV,DLM,LM,V,We}).
\begin{remark}\label{poisonizacion}
{\rm
Let $(\Lambda ,E)$ be a Jacobi structure on a manifold $M$ and
consider on the product manifold $M\times \R$ the 2-vector
$\tilde{\Lambda}$ given by
$$\tilde{\Lambda}=e^{-t}\Big ( \Lambda +\frac{\partial }{\partial
t}\wedge E \Big ) ,$$
where $t$ is the usual coordinate on $\R$. Then, $\tilde{\Lambda}$ defines a
Poisson structure on $M\times \R$ (see \cite{Li2}). The manifold $M\times \R$
endowed with the structure $\tilde{\Lambda}$ is called the {\em Poissonization
of the Jacobi manifold} $(M,\Lambda ,E)$.
}
\end{remark}
\subsection{Lie algebroids}\label{algebroides}
A {\em Lie algebroid} $A$ over a manifold $M$ is a vector bundle
$A$ over $M$ together with a Lie algebra structure on the space
$\Gamma (A)$ of the global cross sections of $A\to M$ and a bundle
map $\rho :A \to TM$, called the {\em anchor map}, such that if we
also denote by $\rho :\Gamma (A) \to \frak X (M)$ the homomorphism of
$C^\infty (M,\R )$-modules induced by the anchor map then:
\begin{enumerate}
\item[(i)] $\rho :(\Gamma (A),\lcf \, ,\, \rcf )\to (\frak X
(M),[\, ,\, ])$ is a Lie algebra homomorphism and
\item[(ii)] for all $f\in C^\infty (M,\R)$ and for all $X ,Y \in
\Gamma (A)$, one has
         $$\lcf X ,fY \rcf=f\lcf X,Y\rcf +(\rho  (X )(f))Y.$$
\end{enumerate}
The triple $(A,\lcf \, ,\, \rcf,\rho )$ is called a {\em Lie algebroid
over} $M$ (see \cite{Mk,Pr}).

A real Lie algebra of finite dimension is a Lie algebroid over a
point. Another trivial example of a Lie algebroid is the triple
$(TM,[\, ,\, ],Id)$, where $M$ is a differentiable manifold and
$Id:TM\to TM$ is the identity map.

If $A$ is a Lie algebroid, the Lie bracket on the sections of $A$ can
be extended to the so-called {\em Schouten bracket} $\lcf \, ,\,
\rcf$ on the space $\Gamma (\wedge ^\ast A)=\oplus_k\Gamma (\wedge
^kA)$ of multi-sections of $A$. The Schouten bracket $\lcf \, ,\,
\rcf :\Gamma (\wedge ^kA)\times \Gamma (\wedge ^{k'}A)\to \Gamma
(\wedge^{k+k'-1}A)$ is characterized by the following conditions:
$\lcf \, ,\, \rcf :\Gamma (A)\times \Gamma (A)\to \Gamma (A)$
coincides with the Lie algebroid bracket, $\lcf X ,f\rcf =\rho (X
)(f)$ for $X \in \Gamma(A)$ and $f\in C^\infty (M,\R)$ and the
properties
$$\lcf P,P' \rcf =(-1)^{kk'}\lcf P',P\rcf ,$$

\vspace{-.9cm}
$$\lcf P,P'\wedge P''\rcf =\lcf P,P'\rcf \wedge P''+(-1)^{k'(k+1)}P'\wedge
\lcf P,P''\rcf ,$$
$$ (-1)^{kk''}\lcf \lcf P,P' \rcf ,P''\rcf +
   (-1)^{k'k''}\lcf \lcf P'',P \rcf ,P'\rcf +
   (-1)^{kk'}\lcf \lcf P',P'' \rcf ,P\rcf =0,$$
holds for all  $P\in \Gamma (\wedge ^kA)$, $P'\in \Gamma (\wedge
^{k'}A)$ and  $P''\in \Gamma (\wedge ^{k''}A)$ (see \cite{MX}).

Next, we will recall the definition of the Lie algebroid cohomology
complex with trivial coefficients. For this
purpose, we recall the definition of the cohomology of a Lie algebra
${\cal A}$ with coefficients in an ${\cal A}$-module (we will follow
\cite{V}).

Let $({\cal A},[\, ,\, ])$ be a real Lie algebra (not necessarily
finite dimensional) and ${\cal M}$ a real vector space endowed with a
$\R$-bilinear multiplication
$${\cal A}\times {\cal M}\to {\cal M},\quad (a,m)\mapsto a\cdot m,$$
such that
$$[a_1,a_2]\cdot m=a_1\cdot (a_2 \cdot m)-a_2\cdot (a_1\cdot m),$$
for $a_1,a_2\in {\cal A}$ and $m\in {\cal M}$. In other words, we
have a representation of ${\cal A}$ on ${\cal M}$. In such a case, a
$k$-linear skew-symmetric mapping $c^k:{\cal A}^k\to {\cal M}$ is
called an  {\em ${\cal M}$-valued k-cochain}. These cochains form a real
vector space $C^k({\cal A};{\cal M})$ and the linear operator $\partial
^k:C^k({\cal A};{\cal M})\to C^{k+1}({\cal A};{\cal M})$ given by
$$
\begin{array}{ccl}
(\partial ^kc^k)(a_0,\ldots ,a_k)&=&\displaystyle \sum
_{i=0}^k(-1)^ia_i\cdot c^k(a_0,\ldots ,\hat{a}_i,\ldots ,a_k)+\\
 & &\displaystyle \sum _{i<j}(-1)^{i+j}c^k([a_i,a_j],a_0,\ldots,
\hat{a}_i,\ldots ,\hat{a}_j,\ldots ,a_k)
\end{array}
$$
defines a coboundary since $\partial ^{k+1}\circ \partial ^k=0$.
Hence we have the corresponding cohomology spaces
$$H^k({\cal A};{\cal M})=\frac{\mbox{ker}\{ \partial ^k:C^k({\cal
A};{\cal M})\to C^{k+1}({\cal A};{\cal M})\}}{\mbox{Im}\{ \partial
^{k-1}:C^{k-1}({\cal A};{\cal M})\to C^k({\cal A};{\cal M})\}}.$$
This cohomology is called {\em the cohomology of the Lie algebra ${\cal A}$
with coefficients in ${\cal M}$, or relative to the given
representation of ${\cal A}$ on ${\cal M}$}.

Now, if $(A,\lcf \, ,\, \rcf ,\rho )$ is a Lie algebroid, we can define
the representation of the Lie algebra $(\Gamma (A),\lcf \,
,\, \rcf )$ on the space $C^\infty (M,$ $\R)$ given by $X \cdot f=\rho
(X )(f)$, for $X \in \Gamma (A)$ and $f\in C^\infty (M,\R )$. We
will denote by $d$ the cohomology operator of the corresponding
cohomology complex. Note that the space of the $k$-cochains which are
$C^\infty (M,\R)$-linear is just $\Gamma (\wedge ^kA^\ast)$, where
$A^\ast$ is the dual bundle to $A$. Moreover, we have that $d(\Gamma
(\wedge ^kA^\ast))\subseteq \Gamma (\wedge ^{k+1}A^\ast)$, for all
$k$, and thus one can consider the subcomplex $(\Gamma (\wedge ^\ast
A^\ast ),d_{|\Gamma (\wedge ^\ast A^\ast )})$. The cohomology
of this subcomplex is the {\em Lie algebroid cohomology with
trivial coefficients} and the restriction of $d$ to $\Gamma (\wedge ^\ast
A^\ast )$ is the {\em differential of the Lie algebroid} $A$ (see
\cite{Mk}).

Using the above definitions, it follows that a 1-cochain $\phi \in
\Gamma (A^\ast)$ is a 1-cocycle if and only if
\begin{equation}\label{condcociclo}
\phi \lcf X ,Y \rcf =\rho (X )(\phi (Y))-\rho (Y)(\phi
(X )), \makebox{ for all }X ,Y \in \Gamma (A).
\end{equation}
To end this section, we will consider two examples of Lie algebroids.

{\bf 1.-} {\em The Lie algebroid} $(TM\times \R ,\makebox{{\bf [}\,
,\, {\bf ]}}, \pi)$

If $M$ is a differentiable manifold, we will exhibit a natural Lie
algebroid structure on the vector bundle $TM\times \R$. First, we
will show some identifications which will be useful in the sequel.

Let $A\to M$ be a vector bundle over $M$. Then, it is clear that
$A\times \R$ is the total space of a vector bundle over $M$.
Moreover, the dual bundle to $A\times \R$ is $A^\ast \times \R$ and the
spaces $\Gamma (\wedge ^r(A\times \R))$ and $\Gamma (\wedge
^k(A^\ast \times \R))$ can be identified with $\Gamma (\wedge^rA)\oplus
\Gamma (\wedge ^{r-1}A)$ and $\Gamma (\wedge ^k A^\ast)\oplus \Gamma
(\wedge ^{k-1}A^\ast)$ in such a way that
\begin{equation}\label{identificaciones}
\begin{array}{ccl}
\kern-10pt(P,Q)((\alpha _1,f_1),\ldots ,(\alpha
_r,f_r))&\kern-10pt=&\kern-10ptP(\alpha _1,\ldots ,\alpha _r)
+\displaystyle \sum _{i=1}^r(-1)^{i+1}f_iQ(\alpha
_1,\ldots,\hat{\alpha} _i,\ldots ,\alpha _r),\\
\kern-10pt(\alpha ,\beta)((X _1,g_1),\ldots ,(X
_k,g_k))&\kern-10pt=&\kern-10pt\alpha (X _1,\ldots
,X _k) +\displaystyle \sum _{i=1}^k(-1)^{i+1}g_i\beta (X
_1,\ldots,\hat{X} _i,\ldots ,X _k),
\end{array}
\end{equation}
for $(P,Q)\in  \Gamma (\wedge^rA)\oplus \Gamma (\wedge ^{r-1}A)$,
$(\alpha ,\beta )\in \Gamma (\wedge ^k A^\ast)\oplus \Gamma
(\wedge ^{k-1}A^\ast)$, $(\alpha _i,f_i)\in \Gamma (A^\ast )\oplus
C^\infty (M,\R )$ and $(X _j,g_j)\in \Gamma (A) \oplus C^\infty (M,\R
)$, with $i\in \{ 1,\cdots ,r\}$ and $j\in \{ 1,\cdots ,k\}$.

Under these identifications, the contractions and the exterior
products are given by
\begin{equation}\label{contraccion}
\begin{array}{lclr}
i_{(\alpha ,\beta )}(P,Q)&=&(i_\alpha P +i_\beta Q,(-1)^ki_\alpha
Q), &\mbox{if }k\leq r,\\
i_{(\alpha ,\beta )}(P,Q)&=&0,& \mbox{if } k>r,\\
i_{(P,Q)}(\alpha ,\beta )&=&(i_P\alpha +i_Q \beta ,(-1)^ri_P\beta ),
&\mbox{if }r\leq k,\\
i_{(P,Q)}(\alpha ,\beta )&=&0,& \mbox{if } r>k,\\
(P,Q)\wedge (P',Q')&=&(P\wedge P',Q\wedge P'+(-1)^rP\wedge Q'),\\
(\alpha ,\beta )\wedge (\alpha ',\beta ')&=&(\alpha \wedge \alpha
',\beta \wedge \alpha '+(-1)^k\alpha \wedge \beta '),
\end{array}
\end{equation}
for $(P',Q')\in \Gamma (\wedge^{r'}A)\oplus \Gamma (\wedge ^{r'-1}A)$ and
$(\alpha ',\beta ')\in \Gamma (\wedge ^{k'} A^\ast)\oplus \Gamma
(\wedge ^{k'-1}A^\ast)$.

Now, suppose that $A$ is the tangent bundle $TM$. Then, the triple
$(A\times \R =TM\times \R ,\makebox{{\bf [}\, ,\, {\bf ]}},
\pi)$ is a Lie algebroid over $M$, where $\pi :TM\times \R \to TM$
is the  canonical projection over the first factor and $\makebox{{\bf
[}\, ,\, {\bf ]}}$ is the bracket given by (see \cite{MMP,V2})
\begin{equation}\label{corchetedeprimerorden}
\makebox{{\bf [}} (X,f),(Y,g)\makebox{{\bf ]}}=([X,Y],X(g)-Y(f)),
\end{equation}
for $(X,f),(Y,g )\in \frak X (M)\times C^\infty (M,\R )\cong \Gamma
(TM\times \R)$. In this case, the dual bundle to $TM\times \R$ is
$T^\ast M\times \R$ and the spaces $\Gamma (\wedge ^r(TM\times \R))$
and $\Gamma (\wedge ^k (T^\ast M\times \R ))$ can be identified with
${\cal V}^r(M)\oplus {\cal V}^{r-1}(M)$ and $\Omega ^k(M)\oplus
\Omega ^{k-1}(M)$. Under these identifications, the differential
$\tilde{\delta}$ of the Lie algebroid is
\begin{equation}\label{unadiferencial}
\tilde{\delta}(\alpha ,\beta )=(\delta \alpha,-\delta \beta )
\end{equation}
and the Schouten bracket $\makebox{{\bf [}}\, ,\, \makebox{{\bf ]}}$
is given by
\begin{equation}\label{unschou}
\makebox{{\bf [}} (P,Q),(P',Q')\makebox{{\bf ]}} =
([P,P'],(-1)^{k+1}[P,Q']-[Q,P']).
\end{equation}

{\bf 2.-} {\em The Lie algebroid $(T^\ast M \kern-1pt \times \R ,\lcf ,\rcf
_{(\Lambda ,E)},\kern-2pt \widetilde{\#}_{(\Lambda ,E)})$ associated with
a Jacobi manifold \kern-1pt $(\kern-2pt M,\Lambda ,E)$}

If $A\to M$ is a vector bundle over $M$ and $P \in
\Gamma (\wedge ^2 A)$ is a 2-section of $A$, we will denote by
$\#_P:\Gamma (A^\ast)\to \Gamma (A)$ the homomorphism of $C^\infty
(M,\R)$-modules given by
\begin{equation}\label{sostenido}
     \beta (\#_P(\alpha ))=P (\alpha ,\beta ),\mbox{ for }\alpha,
                   \beta \in \Gamma (A^\ast ).
\end{equation}
We will also denote by $\#_P:A^\ast\to A$ the corresponding bundle map.

Then, a Jacobi manifold $(M,\Lambda ,E)$ has an associated Lie algebroid
$(T^\ast M \times \R ,\lcf ,\rcf _{(\Lambda ,E)},$ $\widetilde{\#}_{(\Lambda
,E)})$, where $\lcf ,\rcf _{(\Lambda ,E)}$,$\widetilde{\#}_{(\Lambda,E)}$ are
defined by
\begin{equation}\label{ecjacobi}
\begin{array}{cll}
\kern-17.6pt\lcf (\alpha ,f),(\beta ,g)\rcf _{(\Lambda
,E)}&\kern-9pt=&\kern-9pt({\cal L}_{\#_{\Lambda}(\alpha
)}\beta\kern-2pt-\kern-2pt{\cal L}_{\#_{\Lambda}(\beta
)}\alpha\kern-2pt -\kern-2pt \delta (\Lambda
(\alpha ,\beta ))\kern-2pt+\kern-2ptf{\cal
L}_{E}\beta\kern-2pt-\kern-2ptg {\cal L}_{E}\alpha
\kern-2pt-\kern-2pti_{E}(\alpha \wedge \beta),\\
&\kern-9pt &\kern-9pt \Lambda (\beta ,\alpha
)\kern-2pt+\kern-2pt\#_{\Lambda}(\alpha
)(g)\kern-2pt-\kern-2pt\#_{\Lambda}(\beta
)(f)\kern-2pt+\kern-2ptfE(g)\kern-2pt-\kern-2ptg E(f)),\\
& & \\
\kern-17.6pt\widetilde{\#}_{(\Lambda ,E)}(\alpha
,f)&\kern-9pt=&\#_{\Lambda}(\alpha )+fE,
\end{array}
\end{equation}
for $(\alpha ,f),(\beta ,g)\in \Omega ^1(M)\times
C^\infty (M,\R )$, ${\cal L}$ being the Lie derivative operator (see
\cite{KS}). For this algebroid, the differential $d_\ast$ is given by
(see \cite{LLMP,LMP})
\begin{equation}\label{otradiferencial}
d_\ast (P,Q)=(-[ \Lambda ,P ]  +kE\wedge P+\Lambda \wedge Q,[\Lambda
,Q ] -(k-1)E\wedge Q+[E,P]),
\end{equation}
for $(P,Q)\in {\cal V}^k(M)\oplus{\cal V}^{k-1}(M)$.

In the particular case when $(M,\Lambda )$ is a Poisson manifold we
recover, by projection, the Lie algebroid $(T^\ast M,\lcf \, ,\, \rcf
_{\Lambda}, \#_{\Lambda})$, where $\lcf \, ,\, \rcf_{\Lambda}$ is the
bracket of 1-forms defined by (see \cite{BV,CDW,F,V}):
\begin{equation}\label{corchpoisson}
\lcf ,\rcf_{\Lambda}:\Omega ^1(M)\times \Omega ^1(M)\to \Omega ^1(M),\quad
\lcf \alpha ,\beta \rcf _{\Lambda}={\cal
L}_{\#_{\Lambda}(\alpha )}\beta-{\cal L}_{\#_{\Lambda}(\beta )}\alpha
-\delta (\Lambda (\alpha ,\beta )).
\end{equation}
For this algebroid, the differential is the operator $d_\ast
=-[\Lambda ,\cdot \, ]$. This operator was introduced by Lichnerowicz
in \cite{Li1} to define the Poisson cohomology.
\setcounter{equation}{0}
\section{Differential calculus on Lie algebroids in the pre\-sence
of a 1-cocycle}
\subsection{$\phi_0$-differential and $\phi _0$-Lie
derivative}\label{otraseccion}
Let $(A,\lcf ,\rcf ,\rho )$ be a Lie algebroid over $M$ and $\phi _0\in
\Gamma (A^\ast )$ be a 1-cocycle in the Lie algebroid cohomology complex
with trivial coefficients. Using (\ref{condcociclo}), we can define  a
representation $\rho _{\phi _0}:\Gamma (A)\times C^\infty (M,\R )\to
C^\infty (M,\R )$ of the Lie algebra $(\Gamma (A),\lcf \, ,\, \rcf )$
on the space $C^\infty (M,\R )$ given by
\begin{equation}\label{anclaextendida}
      \rho _{\phi _0}(X )f=\rho (X )(f) +\phi _0(X )f,
\end{equation}
for $X \in \Gamma (A)$ and $f\in C^\infty (M,\R )$. Thus, one can consider
the cohomology of the Lie
algebra $(\Gamma (A), \lcf ,\rcf )$ with coefficients in $C^\infty
(M,\R )$ and the subcomplex $\Gamma (\wedge ^\ast A^\ast)$ consisting of the
cochains which are $C^\infty (M,\R)$-linear. The cohomology operator
$d_{\phi _0}:\Gamma (\wedge ^kA^\ast)\to
\Gamma (\wedge ^{k+1}A^\ast)$ of this subcomplex is called the $\phi _0$-{\em
differential} of $A$. We have that
\begin{equation}\label{la1}
d_{\phi _0}\omega =d\omega + \phi _0\wedge \omega,
\end{equation}
where $d$ is the differential of the Lie algebroid $(A,\lcf
\, ,\, \rcf ,\rho)$. As a consequence,
\begin{equation}\label{otramas}
                         d_{\phi _0}1=\phi _0,
\end{equation}
\begin{equation}\label{la2}
d_{\phi _0}(\omega \wedge \omega ')=(d_{\phi _0}\omega )\wedge
\omega '+(-1)^k\omega \wedge (d_{\phi _0}\omega ')-\phi _0\wedge
\omega \wedge \omega ',
\end{equation}
for $\omega \in \Gamma (\wedge ^kA^\ast )$ and $\omega '\in \Gamma
(\wedge ^{k'}A^\ast )$.
\begin{remark}
{\rm
If $\phi _0$ is a closed 1-form on a manifold $M$ then $\phi _0$ is a
1-cocycle for the trivial Lie algebroid $(TM,[\, ,\, ],Id)$ and we
can consider the operator $d_{\phi _0}$. Some results about the
cohomology defined by $d_{\phi _0}$ were obtained in
\cite{GL,LLMP,V0}. These results were used in the study of locally
conformal K\"ahler and locally conformal symplectic structures.
}
\end{remark}
On the other hand, if $k\geq 0$ and $X \in \Gamma (A)$, the
$\phi _0$-{\em Lie derivative} with respect to $X$, $({\cal L}_{\phi _0})_X
:\Gamma (\wedge ^kA^\ast)\to \Gamma (\wedge ^kA^\ast)$, is defined by
\begin{equation}\label{Lieext}
\kern-5pt((({\cal L}_{\phi _0})_X )\omega )(X _1,\ldots ,X
_k)\kern-2pt=\kern-2pt\rho _{\phi
_0}(X)(\omega (X _1,\ldots ,X _k))\kern-1pt-\kern-1pt\sum _{i=1}^k\omega (X
_1,\ldots ,\lcf X , X _i \rcf ,\ldots ,X _k),
\end{equation}
for $\omega \in \Gamma (\wedge ^kA^\ast )$ and $X _1,\ldots ,X _k
\in \Gamma (A)$. It follows that
\begin{equation}\label{rellieext}
({\cal L}_{\phi _0})_X \omega ={\cal L}_X \omega +\phi _0(X
)\omega ,
\end{equation}
${\cal L}$ being the Lie derivative of the Lie algebroid $(A,\lcf
\, ,\,  \rcf ,\rho)$. Thus,
\begin{equation}\label{calcvar}
({\cal L}_{\phi _0})_X =d_{\phi _0}\circ i_X +i_X \circ
d_{\phi _0},
\end{equation}
where $i_X$ is the usual contraction by $X$. Using (\ref{rellieext}) and the
properties of ${\cal L}$  (see \cite{MX}), we deduce that
\begin{proposition}
Let $(A,\lcf ,\rcf ,\rho )$ be a Lie algebroid over $M$ and $\phi _0\in
\Gamma (A^\ast )$ be a 1-cocycle. If $X \in \Gamma (A)$, $f\in
C^\infty(M,\R )$, $\omega \in \Gamma (\wedge ^kA^\ast )$ and $\omega
' \in \Gamma (\wedge ^{k'}A^\ast )$, we have
\begin{equation}\label{derlieprod1}
({\cal L}_{\phi _0})_X (\omega \wedge \omega ')=(({\cal L}_{\phi
_0})_X \omega )\wedge \omega '+\omega \wedge (({\cal L}_{\phi
_0})_X \omega ')-\phi _0(X )\omega \wedge \omega ',
\end{equation}
\begin{equation}
({\cal L}_{\phi _0})_{fX}\omega =f({\cal L}_{\phi _0})_X
\omega+df \wedge i_X \omega .
\end{equation}
\end{proposition}
Now, we will consider the examples of Lie algebroids studied in Section
\ref{algebroides}.

{\bf 1.-} {\em The Lie algebroid} $(TM\times \R ,\makebox{{\bf [}\,
,\, {\bf ]}}, \pi)$

Using (\ref{unadiferencial}), it follows that $\phi _0=(0,1)\in \Omega
^1(M)\times  C^\infty (M,\R )\cong\Gamma (T^\ast M\times \R )$ is a
1-cocycle. Thus, we have the corresponding representation $\pi
_{(0,1)}:(\frak X (M)\times C^\infty (M,\R ))\times C^\infty (M,\R
)\to C^\infty (M,\R )$ of the Lie algebra $(\frak X (M)\times
C^\infty (M,\R),\makebox{{\bf [}\, ,\, {\bf ]}})$ on the space
$C^\infty (M,\R)$ which, in this case, is given by (see
(\ref{anclaextendida}))
\begin{equation}\label{piext}
\pi _{(0,1)}((X,f),g)= X(g)+fg,
\end{equation}
for $(X,f)\in \frak X (M)\times C^\infty (M,\R)$ and $g\in C^\infty (M,\R)$.
From (\ref{identificaciones}), (\ref{unadiferencial}) and
(\ref{la1}), we deduce that the $\phi _0$-differential
$\tilde{\delta}_{\phi _0}=\tilde{\delta} _{(0,1)}$ is given by
\begin{equation}\label{diferencial01}
          \tilde{\delta}_{(0,1)}(\alpha ,\beta )=
              (\delta \alpha ,\alpha -\delta \beta ),
\end{equation}
for $(\alpha ,\beta) \in \Omega ^k(M)\oplus \Omega ^{k-1}(M)$.

\begin{remark}\label{unasnotas}
{\rm
{\em i)} If $(\Lambda ,E)\in \Gamma(\wedge ^2(TM\times \R))$ is a
Jacobi structure on $M$ a long computation, using
(\ref{identificaciones}), (\ref{contraccion}), (\ref{calcvar}) and
(\ref{diferencial01}), shows that the Lie algebroid
bracket $\lcf \, ,\, \rcf _{(\Lambda ,E)}$ and the anchor map
$\widetilde{\#}_{(\Lambda ,E)}$ can be written using the homomorphism
$\# _{(\Lambda ,E)}:\Gamma (T^\ast M \times \R)\to\Gamma (TM \times
\R)$ and the operators ${\cal L}_{(0,1)}$ and
$\tilde{\delta}_{(0,1)}$ as follows
$$
\begin{array}{cll}
\lcf (\alpha ,f),(\beta ,g)\rcf _{(\Lambda ,E)}&=&({\cal
L}_{(0,1)})_{\# _{(\Lambda ,E)}(\alpha ,f)}(\beta ,g)-({\cal
L}_{(0,1)})_{\# _{(\Lambda ,E)}(\beta ,g)}(\alpha ,f)\\
&&-\tilde{\delta}_{(0,1)}\Big ( (\Lambda ,E)((\alpha ,f),(\beta ,g)) \Big )\\
&=&i_{\# _{(\Lambda ,E)}(\alpha ,f)}(\tilde{\delta}_{(0,1)}(\beta ,g))-
i_{\# _{(\Lambda ,E)}(\beta ,g)}(\tilde{\delta}_{(0,1)}(\alpha ,f))\\
&&+\tilde{\delta}_{(0,1)}\Big ( (\Lambda ,E)((\alpha ,f),(\beta ,g))\Big ) ,\\
& & \\
\widetilde{\#}_{(\Lambda ,E)}&=&\pi \circ \#_{(\Lambda,E)}.
\end{array}
$$
Compare equation (\ref{corchpoisson}) with the above expression of
the Lie algebroid bracket $\lcf \, ,\, \rcf _{(\Lambda ,E)}$.

{\em ii)} Let $(A,\lcf \, , \, \rcf ,\rho )$ be a Lie algebroid over $M$ and
$\phi _0\in \Gamma (A^\ast )$ be a 1-cocycle. The homomorphism of
$C^\infty (M,\R)$-modules
\begin{equation}\label{homomorf}
\begin{array}{cccl}
(\rho ,\phi _0):&\Gamma (A)&\to& \frak X(M)\times C^\infty (M,\R)\\
                & X        &\mapsto& (\rho (X),\phi _0 (X)),
\end{array}
\end{equation}
induces a homomorphism between the Lie algebroids
$(A,\lcf \, , \, \rcf ,\rho )$ and $(TM\times \R ,\makebox{{\bf [}\,
,\, {\bf ]}}, \pi)$, that is,
\begin{equation}\label{homomorfismo}
(\rho ,\phi _0)\lcf X,Y\rcf =\makebox{{\bf [}}(\rho ,\phi _0)(X)
,(\rho ,\phi _0)(Y) \makebox{{\bf ]}},\quad
 \pi ((\rho ,\phi _0)(X)) =\rho (X),
\end{equation}
for $X,Y\in \Gamma (A)$.  Moreover, if $(\rho
,\phi _0)^\ast :\Omega ^1(M) \times C^\infty
(M,\R)\to \Gamma (A^\ast)$ is the adjoint homomorphism of $(\rho
,\phi _0)$, then
$(\rho ,\phi _0)^\ast (0,1)=\phi _0.$ As a consequence, the
corresponding homomorphism $(\rho ,\phi _0)^k:\Gamma (\wedge ^kA)\to
\Gamma (\wedge ^k(TM\times \R))\cong
{\cal V}^k(M) \oplus {\cal V}^{k-1}(M)$ and its adjoint homomorphism
$((\rho ,\phi _0)^k)^\ast :\Omega^k(M) \oplus
\Omega^{k-1}(M)\to \Gamma (\wedge ^kA^\ast)$ satisfy
$$((\rho ,\phi _0)^{k+1})^\ast(\tilde{\delta}(\alpha ,\beta))=d(((\rho
,\phi _0)^k)^\ast (\alpha ,\beta )),$$
$$((\rho ,\phi _0)^{k+1})^\ast(\tilde{\delta}_{(0,1)}(\alpha
,\beta))=d_{\phi _0}(((\rho ,\phi _0)^k)^\ast(\alpha ,\beta )).$$
In particular,
\begin{equation}\label{unlema}
(\rho ,\phi _0)^\ast (\tilde{\delta }_{(0,1)}f)=(\rho ,\phi _0)^\ast
(\delta f,f)=d_{\phi _0}f, \mbox{ for }f\in C^\infty (M,\R ).
\end{equation}
}
\end{remark}

{\bf 2.-} {\em The Lie algebroid $(T^\ast M \kern-1pt \times \R ,\lcf ,\rcf
_{(\Lambda ,E)},\kern-2pt \widetilde{\#}_{(\Lambda ,E)})$ associated
with a Jacobi manifold \kern-1pt $(\kern-2pt M,\Lambda ,E)$}

Let $(M,\Lambda ,E)$ be a Jacobi manifold and $(T^\ast M \times \R
,\lcf ,\rcf _{(\Lambda ,E)},\widetilde{\#}_{(\Lambda ,E)})$ the associated
Lie algebroid (see Section \ref{algebroides}). Denote by $d_\ast$ the
differential of $(T^\ast M \times \R ,\lcf ,\rcf _{(\Lambda ,E)},
\widetilde{\#}_{(\Lambda ,E)})$. From (\ref{ecuaciones}) and
(\ref{otradiferencial}), it follows that
$X_0=(-E,0)\in \frak X (M)\times C^\infty (M,\R)\cong\Gamma (TM\times \R
)$ is a 1-cocycle. Using (\ref{otradiferencial}) and (\ref{la1}), we
obtain the following expression for the
$X_0$-differential $d_\ast{}_{X_0}=d_\ast{}_{(-E,0)}$,
\begin{equation}\label{diferencialdual}
\begin{array}{ccl}
d_\ast{}_{(-E,0)}(P,Q)&=&(-[\Lambda ,P]+(k-1)E\wedge P+\Lambda \wedge Q,\\
& & [\Lambda ,Q]-(k-2)E\wedge Q+ [E,P]),
\end{array}
\end{equation}
for $(P,Q)\in {\cal V}^k(M)\oplus{\cal V}^{k-1}(M)$. Note that
$d_\ast{}_{(-E,0)}$ is just the cohomology operator of the {\em
1-differentiable Chevalley-Eilenberg cohomology complex of M} (see
\cite{GL,Li2}).
\subsection{$\phi_0$-Schouten bracket}
In \cite{BV}, a skew-symmetric Schouten bracket was defined for two
multilinear maps of a commutative associative algebra $\frak F$ over
$\R$ with unit as follows. Let ${\cal P}$ and ${\cal P}'$
be skew-symmetric multilinear maps of degree $k$ and $k'$, respectively, and
$f_1,\ldots ,f_{k+k'-1}\in \frak F$. If $A$ is any subset of $\{
1,2,\ldots ,(k+k'-1)\}$, let $A'$ denote its complement and $|A|$ the
number of elements in $A$. If $|A|=l$ and the elements in $A$ are
$\{ i_1,\ldots ,i_l\}$ in increasing order, let us write $f_A$ for the
ordered $k$-uple $(f_{i_1},\ldots ,f_{i_l})$. Furthermore, we write
$\varepsilon _A$ for the sign of the permutation which rearranges the
elements of the ordered $(k+k'-1)$-uple (A',A), in the original order.
Then, the Schouten bracket of ${\cal P}$ and ${\cal P}'$, $\makebox{{\bf [}}
{\cal P} ,{\cal P}' \makebox{{\bf ]}} _{(0,1)}$, is the skew-symmetric
multilinear map of degree $k+k'-1$ given by
$$
\makebox{{\bf [}} {\cal P} ,{\cal P}' \makebox{{\bf ]}} _{(0,1)}
(f_1,\ldots f_{k+k'-1})=\displaystyle \sum_{|A|=k'} \varepsilon _A {\cal P}
({\cal P}' (f_A),f_{A'}) +(-1)^{kk'} \displaystyle \sum_{|B|=k} \varepsilon
_B {\cal P}' ({\cal P} (f_B),f_{B'}).
$$
One can prove that if ${\cal P}$ and ${\cal P}'$ are first-order
differential operators on each of its arguments, so is $\makebox{{\bf
[}} {\cal P} ,{\cal P}' \makebox{{\bf ]}} _{(0,1)}$. In particular, if $M$ is
a differentiable manifold and $\frak F =C^\infty (M,\R )$, we know
that a $k$-linear skew-symmetric first-order differential operator
can be identified with a pair $(P,Q)\in {\cal V}^k(M)\oplus {\cal
V}^{k-1}(M)$ (that is, a $k$-section of $TM\times \R \to M$) in such
a way that
              $$(P,Q)(f_1,\ldots ,f_k)=P(df_1,\ldots
,df_k)+\sum_{i=1}^k(-1)^{i+1}f_i\, Q(df_1,\ldots ,\hat{df_i},\ldots ,df_k),$$
for $f_1,\ldots ,f_k\in C^\infty (M,\R )$. Under the above
identification,  we have that
\begin{equation}\label{unschouext}
\begin{array}{cll}
\makebox{{\bf [}} (P,Q),(P',Q')\makebox{{\bf ]}} _{(0,1)}&=&\Big ( [ P,P']
+(-1)^{k+1}(k-1)P\wedge Q'-(k'-1)Q\wedge P',\\
 & &(-1)^{k+1}[P,Q']-[Q,P']+(-1)^{k+1}(k-k')Q\wedge Q'\Big ) ,
\end{array}
\end{equation}
for $(P,Q)\in {\cal V}^k(M)\oplus{\cal V}^{k-1}
(M)$ and $(P',Q')\in {\cal V}^{k'}(M)\oplus {\cal V}^{k'-1}(M)$. If
$\makebox{{\bf [}} \, , \, \makebox{{\bf ]}}$ is the Schouten bracket
of the Lie algebroid $(TM\times \R ,\makebox{{\bf [}} \, , \,
\makebox{{\bf ]}},\pi )$, an easy computation, using
(\ref{contraccion}), (\ref{unschou}) and (\ref{unschouext}), shows that
\begin{equation}\label{unschouext2}
\begin{array}{ccl}
\kern -10pt \makebox{{\bf [}} (P,Q), (P',Q') \makebox{{\bf
]}}_{(0,1)}&\kern-8pt=&\kern-8pt \makebox{{\bf [}} (P,Q),(P',Q')
\makebox{{\bf ]}} + (-1)^{k+1}(k-1)(P,Q)\wedge (i_{(0,1)}(P',Q'))\\
&&\kern-8pt-(k'-1)(i_{(0,1)}(P,Q))\wedge (P',Q').
\end{array}
\end{equation}
\begin{remark}\label{unasnotas2}
{\rm
{\em i)} Note that $(\Lambda ,E)\in \Gamma (\wedge ^2(TM\times \R))$
defines a Jacobi structure on $M$ if and only if $\makebox{{\bf [}} (\Lambda
,E),(\Lambda ,E) \makebox{{\bf ]}} _{(0,1)}=0$ (see
(\ref{ecuaciones}) and (\ref{unschouext})).

{\em ii)} Using (\ref{diferencialdual}) and (\ref{unschouext}), we have
that the $X_0$-differential $d_\ast{}_{X_0}=d_\ast{}_{(-E,0)}$ of the Lie
algebroid associated with a Jacobi manifold $(M,\Lambda ,E)$ is given by
\begin{equation}\label{diferencialdual2}
d_\ast{}_{(-E,0)}(P,Q)=-\makebox{{\bf [}}(\Lambda ,E),(P,Q)
\makebox{{\bf ]}} _{(0,1)},
\end{equation}
for $(P,Q)\in {\cal V}^k(M)\oplus{\cal V}^{k-1}(M)$. Compare equation
(\ref{diferencialdual2}) with the expression of the differential of
the Lie algebroid associated with a Poisson manifold (see
Section \ref{algebroides}).
}
\end{remark}
Suggested by (\ref{unschouext2}), we prove the following result
\begin{theorem}\label{teoremaschou}
Let $(A,\lcf \, ,\, \rcf ,\rho )$ be a Lie algebroid and $\phi _0\in
\Gamma (A^\ast )$ a 1-cocycle. Then, there exists a unique operation
$\lcf \, ,\, \rcf _{\phi _0}:\Gamma (\wedge ^kA)\times
\Gamma (\wedge ^{k'}A)\to \Gamma (\wedge ^{k+k'-1 }A)$ such that
\begin{equation}\label{anclaext}
\lcf X ,f\rcf _{\phi _0} = \rho _{\phi _0}(X )(f),
\end{equation}
\begin{equation}\label{coinciden}
\lcf X ,Y \rcf _{\phi _0}=\lcf X ,Y \rcf,
\end{equation}
\begin{equation}\label{antisimetria}
\lcf P, P'\rcf _{\phi _0}=(-1)^{kk'}\lcf P', P\rcf _{\phi _0},
\end{equation}
\begin{equation}\label{derivacion}
\lcf P, P'\wedge P''\rcf _{\phi _0}=\lcf P, P'\rcf _{\phi _0}\wedge P''
+(-1)^{k'(k+1)} P'\wedge \lcf P, P''\rcf _{\phi _0}-(i_{\phi
_0}P)\wedge  P'\wedge P'',
\end{equation}
for $f\in C^\infty (M,\R )$, $X ,Y \in \Gamma (A)$, $P \in \Gamma
(\wedge ^k A)$, $P' \in \Gamma (\wedge ^{k'} A)$ and $P'' \in \Gamma
(\wedge ^{k''} A)$.
This operation is given by the general formula
$$
\lcf P, P' \rcf_{\phi _0}=\lcf P, P' \rcf + (-1)^{k+1}(k-1)P\wedge
(i_{\phi _0}P')-(k'-1)(i_{\phi_0}P)\wedge P'.
$$
Furthermore, it satisfies the graded Jacobi identity
\begin{equation}\label{IdJacext}
   (-1)^{kk''}\lcf \lcf P,P' \rcf _{\phi _0},P''\rcf _{\phi _0}+
   (-1)^{k'k''}\lcf \lcf P'',P \rcf _{\phi _0},P'\rcf _{\phi _0}+
   (-1)^{kk'}\lcf \lcf P',P'' \rcf _{\phi _0},P\rcf _{\phi _0}=0.
\end{equation}
\end{theorem}
\prueba We define the operation $\lcf \, ,\, \rcf _{\phi _0}:\Gamma
(\wedge ^kA)\times \Gamma (\wedge ^{k'}A)\to \Gamma (\wedge ^{k+k'-1
}A)$ by
\begin{equation}\label{Schoext}
\lcf P, P' \rcf_{\phi _0}=\lcf P, P' \rcf + (-1)^{k+1}(k-1)P\wedge
(i_{\phi _0}P')-(k'-1)(i_{\phi_0}P)\wedge P'.
\end{equation}
Using (\ref{Schoext}) and properties of the Schouten bracket of
multi-sections of $A$, we
deduce (\ref{anclaext}), (\ref{coinciden}), (\ref{antisimetria}) and
(\ref{derivacion}).

To prove the graded Jacobi identitity, we proceed as follows. If
$\alpha \in \Gamma (A^\ast)$ is a 1-cocycle, we have that
$$i_{\alpha}\lcf X, P'\rcf =\lcf X, i_{\alpha}P'\rcf -i_{d(\alpha (X))}P',$$
for $X\in \Gamma (A)$ and $P'\in \Gamma (\wedge ^{k'} A)$. Using this
relation and the fact that
$$\lcf X_1\wedge \ldots X_k,P'\rcf =\sum_{i=1}^k(-1)^{i+1}X_1\wedge
\ldots \wedge \hat{X}_i\wedge \ldots \wedge X_k\wedge \lcf X_i,P'\rcf ,$$
it follows that
\begin{equation}\label{lema00}
i_{\alpha}\lcf P,P'\rcf=-\lcf i_{\alpha}P,P'\rcf +(-1)^{k+1}\lcf
P,i_{\alpha}P'\rcf ,
\end{equation}
for $P \in \Gamma (\wedge ^k A)$.
From (\ref{Schoext}) and (\ref{lema00}), we deduce that
\begin{equation}\label{lema0}
i_{\phi _0}(\lcf P,P'\rcf_{\phi _0})=-\lcf
i_{\phi _0}P,P'\rcf _{\phi _0}+(-1)^{k+1}\lcf P,i_{\phi _0}P'\rcf_{\phi
_0}.
\end{equation}
On the other hand, we have that
\begin{equation}\label{otraprop}
\lcf P,f\rcf _{\phi _0}=i_{d_{\phi _0}f}P,
\end{equation}
for $f\in C^\infty (M,\R)$. From (\ref{lema00}), (\ref{lema0}) and
(\ref{otraprop}), we obtain that
\begin{equation}\label{IdJacext1}
    \lcf f,\lcf P' ,P''\rcf _{\phi _0}\rcf _{\phi _0}+
    \lcf \lcf f,P' \rcf _{\phi _0},P''\rcf _{\phi _0}+
   (-1)^{k'}\lcf P',\lcf f ,P''\rcf _{\phi _0}\rcf _{\phi _0}=0.
\end{equation}
This proves (\ref{IdJacext}) for $k=0$.

On the other hand if $X\in \Gamma (A)$, using (\ref{lema00}) and the
properties of the Schouten bracket $\lcf \, ,\, \rcf$, it follows that
\begin{equation}\label{IdJacext2}
    \lcf X,\lcf P' ,P''\rcf _{\phi _0}\rcf _{\phi _0}=
    \lcf \lcf X,P' \rcf _{\phi _0},P''\rcf _{\phi _0}+
    \lcf P',\lcf X ,P''\rcf _{\phi _0}\rcf _{\phi _0}.
\end{equation}
We must show that (\ref{IdJacext}) holds, for $k\geq 1$. But, this is
equivalent to prove that (\ref{IdJacext}) holds for $P'\in \Gamma
(\wedge ^{k'}A)$, $P''\in \Gamma (\wedge ^{k''}A)$ and $P=\bar{P}\wedge Y,$
with $\bar{P}\in \Gamma (\wedge ^{k-1}A)$ and $Y\in \Gamma (A)$.

We will proceed by induction on $k$. From (\ref{IdJacext2}), we
deduce that the result is true for $k=1$. Now, assume that
$$
(-1)^{(\bar{k}+1)k''}\lcf \lcf \bar{Q}\wedge Y,P' \rcf _{\phi _0},P''\rcf
_{\phi _0}+ (-1)^{k'k''}\lcf \lcf P'',\bar{Q}\wedge Y \rcf _{\phi
_0},P'\rcf _{\phi _0}+
$$
$$
 (-1)^{(\bar{k}+1)k'}\lcf \lcf P',P'' \rcf _{\phi
_0},\bar{Q}\wedge Y\rcf _{\phi _0}=0,
$$
for $\bar{Q}\in \Gamma (\wedge ^{\bar{k}}A)$, with $\bar{k}\leq k-2$.

Then, we have that
$$
(-1)^{\tilde{k}k''}\lcf \lcf \tilde{Q},P' \rcf _{\phi _0},P''\rcf
_{\phi _0}+ (-1)^{k'k''}\lcf \lcf P'',\tilde{Q} \rcf _{\phi
_0},P'\rcf _{\phi _0}+ (-1)^{\tilde{k}k'}\lcf \lcf P',P'' \rcf _{\phi
_0},\tilde{Q}\rcf _{\phi _0}=0,
$$
for $\tilde{Q}\in \Gamma (\wedge ^{\tilde{k}}A)$, with $\tilde{k}\leq k-1$.

Using this fact, (\ref{lema0}) and (\ref{IdJacext2}), we conclude that
$$
(-1)^{kk''}\lcf \lcf \bar{P}\wedge Y,P' \rcf _{\phi _0},P''\rcf
_{\phi _0}+ (-1)^{k'k''}\lcf \lcf P'',\bar{P}\wedge Y \rcf _{\phi
_0},P'\rcf _{\phi _0}+
$$
$$
 (-1)^{kk'}\lcf \lcf P',P'' \rcf _{\phi
_0},\bar{P}\wedge Y\rcf _{\phi _0}=0.
$$

Finally, if $\lcf \, ,\, \rcf \, \tilde{ }:\Gamma
(\wedge ^kA)\times \Gamma (\wedge ^{k'}A)\to \Gamma (\wedge ^{k+k'-1
}A)$ is an operation which satisfies (\ref{anclaext}),
(\ref{coinciden}), (\ref{antisimetria}) and (\ref{derivacion}), then
it is clear that $\lcf \, ,\, \rcf \, \tilde{ }= \lcf \, ,\, \rcf
_{\phi _0}$.\QED

The operation $\lcf \, ,\, \rcf _{\phi _0}$ is called the {\em $\phi
_0$-Schouten bracket} of $(A,\lcf \, ,\, \rcf ,\rho )$.
Now, if $X \in \Gamma (A)$ and $P\in \Gamma (\wedge ^k A)$, we can
define the {\em $\phi _0$-Lie derivative of $P$ by $X$} as follows
\begin{equation}\label{multlieder}
  ({\cal L}_{\phi _0})_X (P)=\lcf X ,P\rcf _{\phi _0}.
\end{equation}
From Theorem \ref{teoremaschou}, we deduce
\begin{proposition}
Let $(A,\lcf \, ,\, \rcf ,\rho )$ be a Lie algebroid and $\phi _0\in
\Gamma (A^\ast )$ a 1-cocycle. If $f\in C^\infty (M,\R )$, $X \in
\Gamma (A)$, $P \in \Gamma (\wedge ^kA)$ and $P' \in \Gamma (\wedge
^{k'}A)$, we have
\begin{equation}\label{derlieprod2}
({\cal L}_{\phi _0})_X (P\wedge P')=({\cal L}_{\phi _0})_X
(P)\wedge P' + P\wedge ({\cal L}_{\phi _0})_X (P') -\phi _0(X
)P\wedge P',
\end{equation}
\begin{equation}\label{derlieconfunc}
({\cal L}_{\phi _0})_{fX}(P)=f({\cal L}_{\phi _0})_X(P)-X
\wedge i_{df}P.
\end{equation}
\end{proposition}
Using (\ref{Lieext}), (\ref{coinciden}) (\ref{multlieder}) and
(\ref{derlieprod2}), we obtain that
\begin{equation}\label{identificacion}
({\cal L}_{\phi _0})_X(i_{\omega } P)=i_{P}\left ( ({\cal L}_{\phi _0})_X
\omega \right ) + i_{\omega } \left ( ({\cal L}_{\phi _0})_XP\right )
+(k-1)\phi _0(X)i_{\omega }P,
\end{equation}
for $\omega \in \Gamma (\wedge ^kA^\ast)$, $P\in \Gamma
(\wedge  ^kA)$ and  $X \in \Gamma (A)$.

Now, suppose that $(M,\Lambda ,E)$ is a Jacobi manifold. We consider
the 1-jet Lie algebroid $(T^\ast M \times \R ,\lcf \, ,\, \rcf
_{(\Lambda ,E)},\widetilde{\#}_{(\Lambda ,E)})$ associated with the
Jacobi structure $(\Lambda ,E)$ and the 1-cocycle $(-E,0)\in \Gamma
(TM\times \R )\cong\frak X (M)\times C^\infty (M,\R )$. As we know, the
dual bundle $TM\times \R$ admits a Lie algebroid structure
$(\makebox{{\bf [}\, ,\, {\bf ]}},\pi)$ and the pair $(0,1)\in \Gamma
(T^\ast M\times \R)\cong\Omega ^1(M) \times C^\infty (M,\R )$ is a
1-cocycle (see Sections \ref{algebroides} and \ref{otraseccion}).

For the above Lie algebroids and 1-cocycles, we deduce
\begin{proposition}\label{ejemplos}
\begin{itemize}
\item[{\it i)}] If $(X,f) ,(Y,g) \in \Gamma (TM\times \R)\cong\frak X (M)\times
C^\infty (M,\R)$, then
 $$d_\ast{}_{(-E,0)}\makebox{{\bf [}} (X,f) ,(Y,g) \makebox{{\bf ]}}
= \makebox{{\bf [}} (X,f) ,d_\ast{}_{(-E,0)}(Y,g) \makebox{{\bf ]}}
 _{(0,1)}-\makebox{{\bf [}} (Y,g) ,d_\ast{}_{(-E,0)}(X,f)
\makebox{{\bf ]}} _{(0,1)}.$$

\item[{\it ii)}] If ${\cal L}_\ast$ denotes the Lie derivative on the
Lie algebroid $(T^\ast M \times \R ,\lcf ,\rcf _{(\Lambda
,E)},\tilde{\#}_{(\Lambda ,E)})$, then
$$({\cal L}_\ast {}_{(-E,0)})_{(0,1)}(P,Q) +({\cal L}_{(0,1)})_{(-E,0)}(P,Q)
= 0,$$
for $(P,Q)\in \Gamma (\wedge ^k (TM\times \R))\cong{\cal V}^k(M)\oplus
{\cal V}^{k-1}(M)$.
\end{itemize}
\end{proposition}
\prueba {\it i)} It follows from (\ref{diferencialdual2}),
(\ref{antisimetria}) and (\ref{IdJacext}).

{\it ii)} Using (\ref{contraccion}), (\ref{calcvar}),
(\ref{diferencialdual}), (\ref{unschouext}) and (\ref{multlieder}),
we have that
$$
\begin{array}{rcl}
({\cal L}_\ast {}_{(-E,0)})_{(0,1)}(P,Q) +({\cal
L}_{(0,1)})_{(-E,0)}(P,Q)&=&\\
&\kern-300pt =&\kern -150pt d_\ast{}_{(-E,0)}(Q,0)+i_{(0,1)}\Big (
-[\Lambda ,P]+(k-1)E\wedge P+\Lambda \wedge Q,[\Lambda ,Q]\\
&\kern-300pt &\kern-150pt -(k-2)E\wedge Q+[E,P]\Big ) -\Big ( [E,P],[E,Q]
\Big ) =0.
\end{array}
$$

\vspace{-.5cm}
\QED
\begin{remark}
{\rm
{\em i)} Let $(A,\lcf \, ,\, \rcf ,\rho )$ be a Lie algebroid and
$A^\ast$ the dual bundle to $A$. Suppose that $(\lcf \, ,\, \rcf
_\ast,\rho _\ast)$ is a Lie algebroid structure on $A^\ast$. Then,
the pair $(A,A^\ast )$ is said to be a {\em Lie bialgebroid} if
$$d_\ast \lcf X,Y\rcf =\lcf X,d_\ast Y\rcf -\lcf Y,d_\ast X\rcf,$$
for $X,Y\in \Gamma (A)$, where $d_\ast$ is the differential of the
Lie algebroid $(A^\ast ,\lcf \, ,\, \rcf _\ast,\rho _\ast)$ (see
\cite{K-S,MX}).

{\em ii)} Let $(M,\Lambda )$ be a Poisson manifold and $(T^\ast
M,\lcf \, ,\rcf _\Lambda ,\# _\Lambda )$ the associated Lie algebroid
(see Section \ref{algebroides}). If on $TM$ we consider the trivial Lie
algebroid structure, the pair $(TM,T^\ast M)$ is a Lie bialgebroid
(see \cite{MX}).

{\em iii)} If $(M,\Lambda ,E)$ is a Jacobi manifold and on $TM\times
\R$ (resp. $T^\ast M\times \R$) we consider the Lie algebroid
structure $(\makebox{{\bf [}\, ,\, {\bf ]}},\pi)$ (resp. $(\lcf \,
,\, \rcf _{(\Lambda ,E)},\widetilde{\#} _{(\Lambda ,E)})$) then, from
Proposition \ref{ejemplos}, we deduce that the pair $(TM\times
\R,T^\ast M\times \R)$ is not, in general, a Lie bialgebroid (see \cite{V2}).
}
\end{remark}
\setcounter{equation}{0}
\section{Generalized Lie bialgebroids}
\subsection{Generalized Lie bialgebroids and Jacobi structures on the
base space}
Let $A$ be a vector bundle over $M$ and $A^\ast$ the dual bundle
to $A$. Suppose that $(\lcf \, ,\, \rcf ,\rho )$ (resp. $(\lcf \, ,\, \rcf
_\ast,\rho _\ast)$) is a Lie algebroid structure on $A$ (resp. $A^\ast$)
and that $\phi _0\in  \Gamma (A^\ast )$ (resp. $X_0 \in \Gamma (A)$)
is a 1-cocycle in the corresponding Lie algebroid cohomology complex with
trivial coefficients. Then, we will use the following notation:

$\bullet$ $d$ (resp. $d_\ast$) is the differential of $(A,\lcf \, ,\, \rcf
,\rho)$ (resp. $(A^\ast,\lcf \, ,\, \rcf _\ast ,\rho _\ast)$).

$\bullet$ $d_{\phi_0}$ (resp. $d_\ast{}_{X_0}$) is the
$\phi _0$-differential (resp. $X_0$-differential) of $A$ (resp. $A^\ast$).

$\bullet$  ${\cal L}$ (resp. ${\cal L}_\ast$) is the Lie derivative of $A$
(resp. $A^\ast$).

$\bullet$  ${\cal L}_{\phi _0}$ (resp. ${\cal L}_\ast{}_{X_0}$) is the
$\phi _0$-Lie derivative (resp. $X_0$-Lie derivative).

$\bullet$  $\kern-1pt\lcf \, ,\, \rcf _{\phi _0}$ (resp. $\lcf \, ,\, \rcf
_\ast{}_{X _0}$) is the $\phi _0$-Schouten bracket (resp. $\kern-1pt X
_0$-Schouten bracket) on $(\kern-1ptA,\lcf \, ,\, \rcf ,\rho)$ (resp.
$(A^\ast,\lcf \, ,\, \rcf_\ast , \rho _\ast)$).

$\bullet$ $\rho _{\phi _0}:\Gamma (A)\times C^\infty(M,\R)\to C^\infty(M,\R)$
(resp. $\rho _\ast{}_{X_0}:\Gamma (A^\ast)\times C^\infty(M,\R)\to
C^\infty(M,\R)$) is the representation given by (\ref{anclaextendida}).

$\bullet$ $(\rho ,\phi _0):\Gamma (A)\to \frak X (M)\times C^\infty
(M,\R)$ (resp. $(\rho _\ast,X_0):\Gamma (A^\ast)\to \frak X (M)\times C^\infty
(M,\R)$) is the homomorphism of $C^\infty(M,\R)$-modules given by
(\ref{homomorf}) and $(\rho ,\phi _0)^\ast :\Omega ^1 (M)\times C^\infty
(M,\R)\to \Gamma (A^\ast)$ (resp. $(\rho _\ast,X_0)^\ast :\Omega ^1 (M)\times
C^\infty (M,\R)\to \Gamma (A)$) is the adjoint operator of $(\rho
,\phi _0)$ (resp. $(\rho _\ast ,X_0)$).

\begin{definition}\label{genliebialg}
The pair $((A,\phi _0),(A^\ast ,X_0))$ is said to be a generalized
Lie bialgebroid over $M$ if for all $X ,Y \in \Gamma (A)$ and $P\in \Gamma
(\wedge ^kA)$
\begin{equation}\label{condcomp1}
 d_\ast{}_{X_0}\lcf X ,Y \rcf = \lcf X ,d_\ast{}_{X_0}Y \rcf _{\phi
_0}-\lcf Y ,d_\ast{}_{X_0}X \rcf _{\phi _0},
\end{equation}
\begin{equation}\label{condcomp2}
 ({\cal L}_\ast {}_{X_0})_{\phi _0}P+({\cal L}_{\phi _0})_{X_0}P= 0.
\end{equation}
\end{definition}
Using (\ref{anclaextendida}), (\ref{Lieext}), (\ref{coinciden}),
(\ref{multlieder}) and (\ref{derlieprod2}), we obtain that
(\ref{condcomp2}) holds if and only
if
\begin{equation}\label{cond1}
\phi _0 (X_0)=0,\quad \rho (X_0)=-\rho _\ast (\phi _0),
\end{equation}
\begin{equation}\label{cond3}
({\cal L}_\ast)_{\phi _0}X +\lcf X_0 ,X \rcf=0, \mbox{ for
}X \in \Gamma (A).
\end{equation}
Note that (\ref{cond1}) and (\ref{cond3}) follow applying
(\ref{condcomp2}) to $P=f\in C^\infty (M,\R)=\Gamma (\wedge ^0A)$ and
$P=X \in \Gamma (A)$.
\begin{examples}
\begin{itemize}
{\rm
\item[{\it i)}] In the particular case when $\phi _0=0$ and $X_0=0$,
(\ref{condcomp1}) and (\ref{condcomp2}) are equivalent to the condition
$$ d_\ast{}\lcf X ,Y \rcf = \lcf X ,d_\ast Y \rcf -\lcf Y ,d_\ast X \rcf .$$
Thus, the pair $((A,0),(A^\ast,0))$ is a generalized Lie bialgebroid
if and only if the pair $(A,A^\ast)$ is a Lie bialgebroid.
\item[{\it ii)}] Let $(M,\Lambda ,E)$ be a Jacobi manifold. From Proposition
\ref{ejemplos}, we deduce that $\Big ((TM\times \R ,(0,1))$,$(T^\ast M
\times \R ,(-E,0))\Big )$ is a generalized Lie bialgebroid.
}
\end{itemize}
\end{examples}
Next, we will show that if $((A,\phi _0),(A^\ast,X_0))$ is a
generalized Lie bialgebroid over $M$, then $M$ carries an induced
Jacobi structure. First, we will prove some results.
\begin{proposition}\label{relacion}
Let $((A,\phi _0),(A^\ast,X_0))$ be a generalized Lie bialgebroid.
Then,
\begin{equation}\label{4.5}
({\cal L}_\ast{}_{X_0})_{d_{\phi _0}f}X=\lcf X ,d_\ast{}_{X_0}f\rcf ,
\end{equation}
for $X \in \Gamma (A)$ and $f\in C^\infty (M,\R )$.
\end{proposition}
\prueba Using (\ref{la2}) and the derivation law on Lie algebroids,
we obtain that
$$
\begin{array}{ccl}
d_\ast{}_{X_0}\Big ( \lcf X ,fY \rcf \Big )
&=&(d_\ast{}_{X_0}f)\wedge \lcf X ,Y \rcf +fd_\ast{}_{X_0}\lcf X
,Y \rcf -fX_0\wedge \lcf X ,Y \rcf \\
&&+d_\ast{}_{X_0}(\rho (X )(f))\wedge Y +\rho (X )(f)d_\ast{}_{X_0}Y
-\rho (X )(f)X_0\wedge Y,
\end{array}
$$
for $X,Y \in \Gamma (A)$ and $f\in C^\infty (M,\R)$.

On the other hand, from (\ref{la2}),
(\ref{derlieprod2}), (\ref{derlieconfunc}) and (\ref{condcomp1}), we
deduce that
$$
\begin{array}{ccl}
d_\ast{}_{X_0}\Big ( \lcf X ,fY \rcf \Big ) &=&({\cal L}_{\phi
_0})_X (d_\ast{}_{X_0}(fY))-({\cal L}_{\phi _0})_{fY}(d_\ast{}_{X_0}(X)) \\
&=&\Big ( ({\cal L}_{\phi _0})_X (d_\ast{}_{X_0} f)\Big ) \wedge Y
+(d_\ast{}_{X_0}f)\wedge ({\cal L}_{\phi _0})_X Y-\phi_0 (X
)(d_\ast{}_{X_0}f)\wedge Y\\
&& +f({\cal L}_{\phi _0})_X (d_\ast{}_{X_0}Y) +\rho
(X)(f)d_\ast{}_{X_0}Y\\
&& -f\Big (({\cal L}_{\phi _0})_X X_0\wedge Y
+X_0\wedge ({\cal L}_{\phi _0})_X Y -\phi _0(X)X_0\wedge Y
\Big )\\
&&-\rho (X )(f)X_0\wedge Y -f({\cal L}_{\phi _0})_Y
(d_\ast{}_{X_0}X) -i_{df}(d_\ast{}_{X_0}X)\wedge Y.
\end{array}
$$
Thus, using again (\ref{condcomp1}), it follows that
$$
\begin{array}{ccl}
d_\ast{}_{X_0}(\rho (X )(f))\wedge Y&=& \Big ( ({\cal L}_{\phi _0})_X
d_\ast{}_{X_0}f-\phi _0(X )d_\ast{}_{X_0}f-f({\cal L}_{\phi _0})_X X_0\\
&&+f\phi_0(X )X_0-i_{df}(d_\ast{}_{X_0}X )\Big ) \wedge Y,
\end{array}
$$
and so
$$d_\ast{}_{X_0}(\rho (X )(f)) - ({\cal L}_{\phi _0})_X
d_\ast{}_{X_0}f+\phi _0(X )d_\ast{}_{X_0}f+f({\cal L}_{\phi _0})_X X_0-f\phi
_0(X )X_0+i_{df}(d_\ast{}_{X_0}X )=0,$$
which, by (\ref{la2}), (\ref{calcvar}) and (\ref{cond3}), implies
(\ref{4.5}).\QED
\begin{corollary}\label{uncorolario}
Under the same hypothesis as in Proposition \ref{relacion}, we have
\begin{equation}\label{laeccor}
\lcf d_\ast{}_{X_0}g ,d_\ast{}_{X_0}f\rcf =d_\ast{}_{X_0}\Big (
d_{\phi _0}f\cdot d_\ast{}_{X_0}g \Big ) ,
\end{equation}
for all $f,g \in C^\infty (M,\R)$.
\end{corollary}
\prueba
If $\tilde{\delta}_{(0,1)}$ is the operator defined by
(\ref{diferencial01}) then, from (\ref{anclaextendida}), (\ref{Lieext}),
(\ref{calcvar}), (\ref{homomorf}), (\ref{unlema}) and Proposition
\ref{relacion}, we have
$$
\begin{array}{ccl}
\lcf d_\ast{}_{X_0}g,d_\ast{}_{X_0}f\rcf &=&({\cal L}_\ast{}_{X_0})
_{d_{\phi_0}f}(d_\ast{}_{X_0}g) = d_\ast{}_{X_0}\Big ( ({\cal L}_\ast{}_{X_0})
_{d_{\phi _0}f}(g)\Big )\\
&=& d_\ast{}_{X_0}\Big ( \rho _\ast{}_{X_0}(d_{\phi _0}f)(g)\Big ) =
d_\ast{}_{X_0} \Big (
\tilde{\delta}_{(0,1)}g\cdot (\rho _\ast,X_0)(d_{\phi _0}f)\Big )\\
&=&d_\ast{}_{X_0}\Big ( d_{\phi _0}f\cdot (\rho _\ast,X_0)^\ast
(\tilde{\delta}_{(0,1)}g)\Big )
=d_\ast{}_{X_0}\Big ( d_{\phi _0}f\cdot d_\ast{}_{X_0} g\Big ) .
\end{array}
$$
\QED
\begin{remark}\label{notacion}
{\rm Using (\ref{diferencial01}), (\ref{homomorf}) and (\ref{unlema}),
it follows that
$$
\begin{array}{ccl}
d_{\phi _0}f\cdot d_\ast{}_{X_0}g&=&\tilde{\delta}_{(0,1)}f\cdot \Big
( (\rho ,\phi _0)\circ (\rho _\ast, X_0)^\ast \Big )
(\tilde{\delta}_{(0,1)}g)\\
&=&\pi _{(0,1)}\Big ( ((\rho, \phi _0)\circ (\rho _\ast ,X_0)^\ast
)(\tilde{\delta}_{(0,1)}g),f\Big ) ,
\end{array}
$$
where $\pi _{(0,1)}:(\frak X (M)\times C^\infty (M,\R))\times
C^\infty (M,\R)\to C^\infty (M,\R)$ is the representation given by
(\ref{piext}).
}
\end{remark}
Now, we will prove the main result of this section.
\begin{theorem}\label{unteorema}
Let $((A,\phi _0),(A^\ast,X_0))$ be a generalized Lie bialgebroid.
Then, the bracket of functions $\{ \, ,\, \} :C^\infty(M,\R)\times
C^\infty(M,\R)\to C^\infty(M,\R)$ given by
$$\{ f , g \}:=-d_{\phi _0}f\cdot d_\ast{}_{X_0} g,
                \mbox{ for }f,g\in C^\infty (M,\R ),$$
defines a Jacobi structure on $M$.
\end{theorem}
\prueba First of all, we need to prove that $\{ \, , \, \}$ is skew-symmetric.
For that purpose, we first show that
\begin{equation}\label{condicionantisim}
\tilde{\delta}_{(0,1)}f\cdot \Big ( ((\rho ,\phi _0)\circ (\rho _\ast
,X_0) ^\ast)(0,1)\Big )=
-\tilde{\delta}_{(0,1)}1\cdot \Big ( ((\rho ,\phi _0)\circ (\rho _\ast
,X_0) ^\ast)(\delta f, 0)\Big ),
\end{equation}
for $f\in C^\infty (M,\R )$. We have that (see (\ref{diferencial01}),
(\ref{homomorf}) and (\ref{cond1}))
$$\tilde{\delta}_{(0,1)}f\cdot \Big ( ((\rho ,\phi _0)\circ (\rho _\ast
,X_0) ^\ast)(0,1)\Big ) = (\delta f,f)\cdot (\rho (X_0),0)
= \rho (X _0)(f)=-\rho _\ast(\phi _0)(f).$$
On the other hand, from (\ref{homomorf}) and (\ref{cond1}), we deduce that
$$
-\tilde{\delta}_{(0,1)}1\cdot \Big ( ((\rho ,\phi _0)\circ (\rho _\ast
,X_0) ^\ast)(\delta f, 0)\Big )=
- \phi _0\cdot (\rho _\ast ,X_0)^\ast (\delta f,0)
)=-\rho _\ast (\phi _0)(f).
$$
Thus, we deduce (\ref{condicionantisim}). Using (\ref{diferencial01}),
(\ref{condicionantisim}) and Remark \ref{notacion}, we obtain that
\begin{equation}\label{segundacond}
\begin{array}{ccl}
\{ f,f\}&=&-\Big ( ((\rho ,\phi_0)\circ (\rho _\ast ,X_0)^\ast )
(\delta f,0)\Big ) \cdot (\delta f,0).
\end{array}
\end{equation}
Now, we will prove that
         $$\Big ( ((\rho ,\phi_0) \circ (\rho _\ast ,X_0)^\ast)(\delta
                      f,0)\Big ) \cdot (\delta f,0)=0.$$
From (\ref{la2}), (\ref{segundacond}), Corollary \ref{uncorolario} and Remark
\ref{notacion}, it follows that
$d_\ast{}_{X_0}\Big ( ((\rho ,\phi_0) \circ (\rho _\ast
,X_0)^\ast )$ $(\delta f^2,0)\cdot (\delta f^2,0) \Big )=0$. Then,
\begin{equation}\label{unaigualdad}
\begin{array}{ccl}
0&=&\Big ( ((\rho ,\phi_0) \circ (\rho _\ast ,X_0)^\ast )(\delta
f,0)\cdot (\delta f,0)\Big ) (d_\ast{}_{X_0}f^2 -f^2X_0)\\
&=&\Big (  ((\rho ,\phi_0) \circ (\rho _\ast ,X_0)^\ast )(\delta
f,0)\cdot (\delta f,0) \Big ) d_\ast f^2\\
&=&2f\Big ( ((\rho ,\phi_0) \circ (\rho _\ast ,X_0)^\ast)(\delta
f,0)\cdot (\delta f,0)\Big ) d_\ast f.
\end{array}
\end{equation}
On the other hand, in general, $d_\ast g=(\rho _\ast ,X _0)^\ast
(\delta g,0).$ Thus, using (\ref{unaigualdad}),
   $$f\Big ( ((\rho ,\phi_0) \circ (\rho _\ast ,X_0)^\ast )(\delta
      f,0)\cdot (\delta f,0)\Big ) ^2 =0, \quad \mbox{for all }f.$$
This implies that
       $$\Big ( ((\rho ,\phi_0) \circ (\rho _\ast ,X_0)^\ast )(\delta
               f,0)\Big ) \cdot (\delta f,0)=0,$$
as we wanted to prove. Therefore, we conclude that $\{ f ,f\} =0$,
for all $f\in C^\infty (M,\R)$, that is, $\{ \, ,\, \}$ is
skew-symmetric.

From (\ref{otramas}) and (\ref{la2}), we deduce that $\{ \, , \}$ is a
first-order differential operator on each of its arguments.

Now, let us prove the Jacobi identity. Using (\ref{laeccor}), we have
$$d_{\phi _0}h\cdot \lcf d_\ast{}_{X_0}g,d_\ast{}_{X_0}f\rcf
=-d_{\phi _0}h\cdot d_\ast{}_{X_0}(\{ f,g\}).$$
Thus, from (\ref{homomorfismo}) and (\ref{unlema}), we deduce that
$$\tilde{\delta}_{(0,1)}h\cdot \makebox{{\bf [}} (\rho ,\phi
_0)((\rho _\ast ,X_0)^\ast (\tilde{\delta}_{(0,1)}g)),(\rho ,\phi _0)
((\rho _\ast,X_0)^\ast (\tilde{\delta}_{(0,1)}f)) \makebox{{\bf ]}}=
-d_{\phi _0}h\cdot d_\ast{}_{X_0}(\{ f,g\} ),$$
or equivalently
$$\pi _{(0,1)}\Big ( \makebox{{\bf [}}((\rho ,\phi _0)\circ(\rho _\ast
,X_0)^\ast )(\tilde{\delta}_{(0,1)}(g)),((\rho ,\phi _0)\circ(\rho _\ast
,X_0)^\ast )(\tilde{\delta}_{(0,1)}(f))\makebox{{\bf ]}},h\Big
)\kern-.7pt =\kern-.7pt-d_{\phi _0}h\cdot d_\ast{}_{X_0}(\{ f,g\} ). $$
Consequently, since $\pi _{(0,1)}$ is a representation of the Lie algebra
$(\frak X (M)\times C^\infty (M,\R),\makebox{{\bf [}\, ,\, {\bf ]}})$
on the space $C^\infty (M,\R)$, this implies that (see Remark \ref{notacion})
$$\{ f,\{ g,h\} \} +\{ g,\{ h,f\} \} +\{ h,\{ f,g\} \} =0.$$
\QED

From (\ref{la1}) and (\ref{cond1}), we have that
\begin{equation}\label{desarrollocorch}
\{ f,g\} =-df\cdot d_\ast g +f\rho (X_0)(g)-g\rho (X_0)(f),
\end{equation}
for $f,g\in C^\infty (M,\R)$. Since the differential $d$ is a
derivation with respect to the ordinary multiplication of functions
we have that the map $(f,g)\mapsto df\cdot d_\ast g$, for $f,g\in
C^\infty(M,\R)$, is also a derivation on each of its arguments. Thus,
we can define the 2-vector $\Lambda \in {\cal V}^2(M)$
characterized by the relation
\begin{equation}\label{Lambda}
\Lambda (\delta f,\delta g)=-df \cdot d_\ast g=dg\cdot d_\ast f,
\end{equation}
for $f,g\in C^\infty (M,\R)$, and the vector field $E\in \frak X(M)$ by
\begin{equation}\label{E}
E=\rho (X_0)=-\rho _\ast (\phi _0).
\end{equation}
From (\ref{desarrollocorch}), we obtain that
   $$\{f,g\} =\Lambda (\delta f,\delta g)+fE(g)-gE(f),$$
for $f,g\in C^\infty (M,\R)$. Therefore, the pair $(\Lambda ,E)$
is the Jacobi structure induced by the Jacobi bracket $\{ \, ,\, \}$.

If $(A,A^\ast)$ is a Lie bialgebroid, then the pair $((A,0),(A^\ast,0))$ is
a generalized Lie bialgebroid and, by Theorem \ref{unteorema}, a
Jacobi structure $(\Lambda, E)$ can be defined on the base space $M$.
Since $\phi _0=X_0=0$, we deduce that $E=0$, that is, the Jacobi
structure is Poisson, which implies a well known result (see
\cite{MX}): given a Lie bialgebroid $(A,A^\ast)$ over $M$, the base space $M$
carries an induced Poisson structure.
\subsection{Lie bialgebroids and generalized Lie bialgebroids}
Along this section we will use the same notation as in the precedent one.
\subsubsection{Time-dependent sections of a Lie algebroid}
Let $(A,\lcf \, , \, \rcf, \rho )$ be a Lie algebroid over $M$. Then,
the space of sections $\Gamma (\tilde{A})$ of the vector bundle
$\tilde{A}=A\times \R \to M\times \R$ can be identified with the set
of time-dependent sections of $A\to M$. Under this identification,
the Lie bracket $\lcf \, , \, \rcf$ induces, in a natural way, a Lie
bracket on $\Gamma (\tilde{A})$ which is also denoted by $\lcf \, ,
\, \rcf$. In fact, if $\tilde{X},\tilde{Y}\in \Gamma
(\tilde{A})$ then $\lcf \tilde{X} ,\tilde{Y} \rcf$ is the
time-dependent section of $A\to M$ given by
$$\lcf \tilde{X}, \tilde{Y}\rcf (x,t)=\lcf
\tilde{X}_t,\tilde{Y}_t\rcf (x),$$
for all $(x,t)\in M\times \R$, where $\tilde{X}_t$ and $\tilde{Y}_t$
are the sections of $A\to M$ defined by
$$\tilde{X}_t(y)=\tilde{X}(y,t),\quad \tilde{Y}_t(y)=\tilde{Y}(y,t),$$
for all $y\in M$. The anchor map $\rho :A\to TM$ induces a bundle map
from $\tilde{A}$
into $TM  \subseteq T(M\times \R)\cong TM\oplus T\R$ which we also
denote by $\rho$. A direct computation shows that the triple
$(\tilde{A},\lcf \, , \,
\rcf, \rho )$ is a Lie algebroid over $M\times \R$.

Now, denote by $d$ the differentials of the Lie algebroids $(A,\lcf
\, , \, \rcf, \rho )$ and $(\tilde{A},\lcf \, , \, \rcf, \rho )$. Then, if
$\tilde{\omega} \in \Gamma (\wedge ^k\tilde{A}^\ast)$ and $(x,t)\in
M\times \R$, $d\tilde{\omega} \in
\Gamma (\wedge ^{k+1}\tilde{A}^\ast)$ and
            $$(d\tilde{\omega})(x,t)=(d\tilde{\omega} _t)(x).$$
We also denote by $\lcf \, , \, \rcf$ the Schouten bracket of the Lie
algebroid $(\tilde{A},\lcf \, , \, \rcf, \rho )$.

On the other hand, for $\tilde{P}\in \Gamma (\wedge ^k \tilde{A})$ or
$\tilde{\omega}\in \Gamma (\wedge ^k \tilde{A}^\ast)$, one can define
its derivative with respect to the time
$$\displaystyle \frac{\partial\tilde{P}}{\partial t}\in \Gamma
(\wedge ^k\tilde{A}),\qquad
 \displaystyle \frac{\partial\tilde{\omega}}{\partial t}\in \Gamma
(\wedge ^k\tilde{A}^\ast ).$$
Thus, we have two $\R$-linear operators of degree zero
$$\frac{\partial}{\partial t}:\Gamma (\wedge ^k\tilde{A})\to \Gamma
(\wedge ^k\tilde{A}),\qquad  \frac{\partial}{\partial t}:\Gamma
(\wedge ^k\tilde{A}^\ast )\to \Gamma (\wedge ^k\tilde{A}^\ast ), $$
which have the following properties
\begin{equation}\label{tercerlema0}
\frac{\partial}{\partial t}(\tilde{P}\wedge \tilde{Q})=\frac{\partial \tilde{P}}{\partial
t}\wedge \tilde{Q} +\tilde{P}\wedge \frac{\partial \tilde{Q}}{\partial t},\quad
\frac{\partial}{\partial t}(\tilde{\omega} \wedge \tilde{\mu} )=\frac{\partial
\tilde{\omega}}{\partial t}\wedge \tilde{\mu} +\tilde{\omega} \wedge
\frac{\partial \tilde{\mu}}{\partial t},
\end{equation}
\begin{equation}\label{tercerlema1}
\frac{\partial}{\partial t}\lcf \tilde{P},\tilde{Q} \rcf  =\lcf
\frac{\partial \tilde{P}}{\partial t},\tilde{Q} \rcf +\lcf
\tilde{P},\frac{\partial \tilde{Q}}{\partial t}\rcf ,
\end{equation}
\begin{equation}\label{tercerlema2}
d\left ( \frac{\partial \tilde{\omega}}{\partial t} \right
)=\frac{\partial}{\partial t}\left ( d\tilde{\omega}\right ),
\end{equation}
for $\tilde{P}\in \Gamma (\wedge ^k \tilde{A})$, $\tilde{Q}\in \Gamma
(\wedge ^{r} \tilde{A})$, $\tilde{\omega} \in \Gamma (\wedge
^k\tilde{A}^\ast)$ and $\tilde{\mu} \in \Gamma (\wedge ^r\tilde{A}^\ast)$.

Next, suppose that $\phi _0\in \Gamma (A^\ast )$ is a 1-cocycle.
In \cite{IM} we proved that the vector bundle $\tilde{A} \to M\times
\R$ admits a Lie algebroid structure $(\lcf \, ,\, \rcf \, \hat{ }\, ^{\phi
_0},\hat{\rho}^{\phi _0})$
where, under the above identifications, $\lcf \, , \, \rcf \, \hat{ }\,
^{\phi _0}$ and $\hat{\rho}^{\phi _0}$ are defined by
\begin{equation}\label{corchtilde}
\begin{array}{c}
\lcf \tilde{X},\tilde{Y}\rcf \, \hat{ }\, ^{\phi _0}=e^{-t}\Big( \lcf \tilde{X},
\tilde{Y}\rcf + \phi _0(\tilde{X})(\frac{\partial \tilde{Y}}{\partial
t}-\tilde{Y})-\phi _0(\tilde{Y})(\frac{\partial \tilde{X}}{\partial
t}-\tilde{X})\Big ) ,\\
\\
\hat{\rho}^{\phi _0}(\tilde{X})= e^{-t}\Big( \rho (\tilde{X})+\phi_0(\tilde{X})
\frac{\partial}{\partial t}\Big ) ,
\end{array}
\end{equation}
for all $\tilde{X},\tilde{Y}\in \Gamma(\tilde{A})$. If $\hat{d}^{\phi _0}$ is
the differential of the Lie algebroid $(\tilde{A},\lcf ,\rcf \, \hat{
}\, ^{\phi _0},\hat{\rho}^{\phi _0})$, we have
\begin{equation}\label{difcomp}
\hat{d}^{\phi _0}\tilde{f}=e^{-t}\Big (d\tilde{f}+\frac{\partial \tilde{f}}{\partial t}
\phi _0\Big ),
\end{equation}
\begin{equation}\label{primerlema}
\hat{d}^{\phi _0} \tilde{\phi}=e^{-t}\left ( d_{\phi _0}\tilde{\phi} +\phi _0\wedge
\frac{\partial \tilde{\phi}}{\partial t}\right ),
\end{equation}
for $\tilde{f}\in C^\infty (M\times \R,\R)$ and $\tilde{\phi }\in
\Gamma(\tilde{A}^\ast)=\Gamma (A^\ast\times \R )$.

In a similar way, we can define the bracket $\lcf \,,\, \rcf \, \bar{
}\, ^{\phi _0}$
on $\Gamma (\tilde{A})$ and the map $\bar{\rho}^{\phi _0}:\Gamma (\tilde{A})
\to \frak X (M\times \R)$ by
\begin{equation}\label{corchbarra}
\begin{array}{c}
\lcf \tilde{X},\tilde{Y}\rcf \, \bar{ }\, ^{\phi _0}= \lcf \tilde{X},
\tilde{Y}\rcf + \phi _0(\tilde{X})\frac{\partial \tilde{Y}}{\partial t}
-\phi _0(\tilde{Y})\frac{\partial \tilde{X}}{\partial t},\\
\\
\bar{\rho}^{\phi _0}(\tilde{X})= \rho (\tilde{X})+\phi _0(\tilde{X})
\frac{\partial}{\partial t},
\end{array}
\end{equation}
for all $\tilde{X},\tilde{Y}\in \Gamma(\tilde{A})$. Moreover, the
bundle map $\Psi :\tilde{A}\to \tilde{A},(v,t)\mapsto (e^tv,t)$,
is an isomorphism of vector bundles and
$$\hat{\rho}^{\phi _0}\circ \Psi =\bar{\rho}^{\phi _0},\qquad
\Psi \lcf \tilde{X},\tilde{Y}\rcf \, \bar{ }\, ^{\phi _0}=\lcf \Psi
\tilde{X}, \Psi \tilde{Y}\rcf \, \hat{ }\, ^{\phi _0}.$$
Thus,
\begin{proposition}\label{isomorf}
Let $(A,\lcf \, , \, \rcf, \rho )$ be a Lie algebroid and $\phi _0\in
\Gamma (A^\ast )$ a 1-cocycle. Then:
\begin{itemize}
\item[{\it i)}] The triples $(\tilde{A},\lcf \, ,\, \rcf \, \hat{ }\, ^{\phi
_0},$ $\hat{\rho}^{\phi _0})$ and $(\tilde{A},\lcf \, ,\, \rcf \,
\bar{ }\, ^{\phi _0}, \bar{\rho}^{\phi _0})$ are Lie algebroids over
$M\times \R$.
\item[{\it ii)}] The map $\Psi :\tilde{A}\to \tilde{A}$ defines an isomorphism
between the Lie algebroids $(\tilde{A},\lcf \, ,\, \rcf \, \bar{
}\, ^{\phi _0},$ $\bar{\rho}^{\phi _0})$ and $(\tilde{A},\lcf \, ,\,
\rcf \, \hat{ }\, ^{\phi _0},\hat{\rho}^{\phi _0})$.
\end{itemize}
\end{proposition}

If $\bar{d}^{\phi _0}$ is the differential of the algebroid $(\tilde{A},\lcf \,
,\, \rcf \, \bar{ }\, ^{\phi _0},\bar{\rho}^{\phi _0})$, we have
\begin{equation}\label{difcomp2}
\bar{d}^{\phi _0}\tilde{f}=d\tilde{f}+\frac{\partial \tilde{f}}{\partial t}\phi _0,
\end{equation}
\begin{equation}\label{difcomp22}
\bar{d}^{\phi _0} \tilde{\phi}= d\tilde{\phi} +\phi _0\wedge
\frac{\partial \tilde{\phi}}{\partial t},
\end{equation}
for $\tilde{f}\in C^\infty (M\times \R,\R)$ and $\tilde{\phi }\in
\Gamma(\tilde{A}^\ast)=\Gamma (A^\ast\times \R)$. Moreover,
\begin{equation}\label{segundolema}
\lcf \tilde{X},\tilde{P} \rcf \, \bar{ }\, ^{\phi _0}=\lcf \tilde{X},
\tilde{P}\rcf _{\phi _0}+\phi _0(\tilde{X})\left ( \tilde{P}+\frac{\partial
\tilde{P}}{\partial t} \right )-\frac{\partial \tilde{X}}{\partial
t}\wedge i_{\phi _0}\tilde{P},
\end{equation}
for $\tilde{X}\in \Gamma(\tilde{A})$ and $\tilde{P}\in \Gamma (\wedge
^2\tilde{A})$.

Now, let $A\to M$ be a vector bundle over a manifold $M$ and suppose
that $\lcf \, ,\, \rcf :\Gamma (A)\times \Gamma (A) \to \Gamma (A)$
is a bracket on the space $\Gamma (A)$, that $\rho :\Gamma (A)\to
\frak X (M)$ is a homomorphism of $C^\infty (M,\R )$-modules and that
$\phi _0$ is a section of the dual bundle $A^\ast$.

We can define the bracket $\lcf \, ,\, \rcf \, \bar{ }\, ^{\phi _0}:\Gamma
(\tilde{A})\times \Gamma (\tilde{A})\to \Gamma (\tilde{A})$ on the
space $\Gamma (\tilde{A})$ and the homomorphism  of $C^\infty (M\times
\R ,\R )$-modules $\bar{\rho}^{\phi _0}: \Gamma (\tilde{A})\to \frak X (M\times
\R)$ given by (\ref{corchbarra}).
\begin{proposition}\label{otrapropo}
If the triple $(\tilde{A},\lcf \, ,\, \rcf \, \bar{ }\, ^{\phi
_0},\bar{\rho}^{\phi _0})$ is a Lie
algebroid on $M\times \R$ then the triple $(A,\lcf \, , \, \rcf,
\rho )$ is a Lie algebroid on $M$ and the section $\phi _0$ is a 1-cocycle.
\end{proposition}
\prueba From (\ref{corchbarra}), it follows that
          $\lcf X,Y\rcf \, \bar{ }\, ^{\phi _0}=\lcf X,Y\rcf ,$
for $X,Y\in \Gamma (A)$. Thus, we have that the bracket $\lcf \, ,\,
\rcf$ defines a Lie algebra structure on $\Gamma (A)$.

On the other hand, since $\bar{\rho}^{\phi _0}\lcf X,Y\rcf \, \bar{
}\, ^{\phi _0}=[\bar{\rho}^{\phi _0}(X),\bar{\rho}^{\phi _0}(Y)],$
we deduce that (see (\ref{corchbarra}))
$$
\rho\lcf X,Y\rcf =[\rho(X),\rho(Y)],\quad
\phi _0\lcf X,Y\rcf =\rho (X)(\phi _0(Y))-\rho (Y)(\phi _0(X)).
$$
Finally, if $f\in C^\infty (M,\R )$ then, using (\ref{corchbarra})
and the fact that $\lcf X,fY\rcf \, \bar{ }\, ^{\phi _0}=f\lcf X,
Y\rcf \, \bar{ }\, ^{\phi _0} + (\bar{\rho}^{\phi _0}(X)(f))Y$, we obtain that
      $$\lcf X,fY\rcf =f\lcf X, Y\rcf + (\rho(X)(f))Y.$$
\QED

From Propositions \ref{isomorf} and \ref{otrapropo}, we conclude
\begin{proposition}\label{equivalencia}
Let $A\to M$ be a vector bundle over a manifold $M$. Suppose
that $\lcf \, ,\, \rcf :\Gamma (A)\times \Gamma (A) \to \Gamma (A)$
is a bracket on the space $\Gamma (A)$, that $\rho :\Gamma (A)\to
\frak X (M)$ is a homomorphism of $C^\infty (M,\R )$-modules and that
$\phi _0$ is a section of the dual bundle $A^\ast$. If $\lcf \, ,\,
\rcf \, \hat{ }\, ^{\phi _0}:\Gamma (\tilde{A})\times \Gamma
(\tilde{A})\to \Gamma (\tilde{A})$ and $\hat{\rho}^{\phi _0}: \Gamma
(\tilde{A})\to \frak X (M\times \R)$ (respectively, $\lcf \, ,\, \rcf
\, \bar{ }\, ^{\phi _0}:\Gamma (\tilde{A})\times \Gamma
(\tilde{A})\to \Gamma (\tilde{A})$ and
$\bar{\rho}^{\phi _0}: \Gamma (\tilde{A})\to \frak X (M\times
\R)$) are the bracket on $\Gamma (\tilde{A})$ and the homomorphism of
$C^\infty (M\times \R ,\R )$-modules given by (\ref{corchtilde})
(respectively, (\ref{corchbarra})) then the following conditions are
equivalent:

i) The triple $(A,\lcf \, , \, \rcf, \rho )$ is a Lie algebroid and
$\phi _0$ is a 1-cocycle.

ii) The triple $(\tilde{A},\lcf \, ,\, \rcf \, \hat{ }\,
^{\phi _0},\hat{\rho}^{\phi _0})$ is a Lie algebroid.

iii) The triple $(\tilde{A},\lcf \, ,\, \rcf \, \bar{ }\,
^{\phi _0},\bar{\rho}^{\phi _0})$
is a Lie algebroid.
\end{proposition}
\subsubsection{Lie bialgebroids and generalized Lie bialgebroids}
First of all, we will prove a general result which will be useful in
the sequel.

Suppose that $(A_i,\lcf \, ,\, \rcf _i,\rho _i )$, $i=1,2$, are two
Lie algebroids over $M$ such that the dual bundles $A_1^\ast$ and
$A_2^\ast$ are Lie algebroids with Lie algebroid structures
$(\lcf \, ,\, \rcf _1{}_\ast,\rho _1{}_\ast)$ and $(\lcf \, ,\,
\rcf _2{}_\ast,\rho _2{}_\ast)$, respectively.
\begin{proposition}\label{unaprop}
Let $\Phi :A_1\to A_2$ be a Lie algebroid isomorphism
such that its adjoint homomorphism $\Phi ^\ast :A_2^\ast \to A_1^\ast$ is
also a Lie algebroid isomorphism. Then, if
$(A_1,A_1^\ast )$ is a Lie bialgebroid, so is $(A_2,A_2^\ast )$.
\end{proposition}
\prueba Denote also by $\Phi: \wedge ^kA_1\to \wedge ^kA_2$ the
isomorphism between the vector bundles $\wedge ^kA_1\to M$ and
$\wedge ^kA_2\to M$ induced by $\Phi :A_1\to A_2$. If $\Phi :\Gamma (\wedge
^kA_1)\to \Gamma (\wedge ^kA_2)$ is the
corresponding isomorphism of $C^\infty (M,\R )$-modules, we have that
$$\Phi (P)(\psi _1,\ldots, \psi _k)=P(\Phi ^\ast \psi _1,\ldots ,\Phi
^\ast \psi _k),$$
$$\Phi (X_1\wedge \ldots \wedge X_k)=\Phi (X_1)\wedge \ldots \wedge
\Phi (X_k),$$
for $P\in \Gamma (\wedge ^kA_1), \psi _1,\ldots ,\psi _k\in \Gamma
(A_2^\ast )$ and $X_1,\ldots ,X_k\in \Gamma (A_1)$. Thus, using that
$\Phi$ and $\Phi ^\ast$ are Lie algebroid isomorphisms, it follows that
\begin{equation}\label{4.24'}
        d_2{}_\ast (\Phi (X_1))=\Phi (d_1{}_\ast X_1),\qquad
    \Phi \lcf X_1,P_1\rcf _1=\lcf \Phi (X_1),\Phi (P_1)\rcf _2,
\end{equation}
for $X_1\in \Gamma (A_1)$ and $P_1\in \Gamma (\wedge ^k A_1)$, where $d_1{}
_\ast$ (resp. $d_2{}_\ast$) is the differential of $(A_1^\ast ,\lcf \, ,\,
\rcf _1{}_\ast,$ $\rho _1{}_\ast )$ (resp. $(A_2^\ast ,\lcf \, ,\, \rcf _2{}
_\ast ,\rho _2{}_\ast)$).

Now, if $X_2,Y_2\in \Gamma (A_2)$ then there exist $X_1,Y_1\in \Gamma
(A_1)$ such that $Y_i=\Phi (X_i)$, for $i=1,2$. Therefore, from
(\ref{4.24'}) and since $(A_1,A^\ast _1)$ is a Lie bialgebroid, we
obtain that
$$\begin{array}{ccl}
d_2{}_\ast\lcf X_2,Y_2\rcf _2&=& (d_2{}_\ast \circ \Phi )\lcf X_1,Y_1\rcf _1
 =\Phi (\lcf X_1, d_1{}_\ast Y_1\rcf _1-\lcf Y_1,
d_1{}_\ast X_1\rcf _1)\\
&=&\lcf X_2, d_2{}_\ast Y_2\rcf _2-\lcf Y_2, d_2{}_\ast X_2 \rcf _2.
\end{array}
$$
Consequently, $(A_2,A^\ast_2)$ is a Lie bialgebroid.\QED

Next, assume that $(M,\Lambda ,E)$ is a Jacobi manifold. Consider on
$A=TM\times \R$
and on $A^\ast =T^\ast M\times \R$ the Lie algebroid structures
$(\makebox{{\bf [}\, ,\, {\bf ]}}, \pi)$ and $(\lcf ,\rcf _{(\Lambda
,E)},\widetilde{\#}_{(\Lambda ,E)})$, respectively. Then, the pair
$((A,\phi _0=(0,1)),(A^\ast ,X_0=(-E,0)))$ is a generalized Lie bialgebroid.

On the other hand, the map $\Phi :\tilde{A}=A\times \R \to T(M\times
\R)$ defined by
      $$\Phi ((v_{x_0},\lambda _0),t_0)= v_{x_0}+\lambda
           _0\frac{\partial}{\partial t}_{|t_0},$$
for $x_0\in M$, $v_{x_0}\in T_{x_0}M$ and $\lambda _0,t_0\in \R$,
induces an isomorphism between the vector bundles $A\times \R\to
M\times \R$ and $T(M\times \R)\to M\times \R$. Thus,
$\tilde{A}=A\times \R$ can be identified with $T(M\times \R)$ and,
under this identification, the Lie algebroid structure
$(\makebox{{\bf [}\, ,\, {\bf ]}}\, \bar{ }\, ^{\phi _0},
\bar{\pi}^{\phi _0})$ is just the
trivial Lie algebroid structure on $T(M\times \R)$ (see
(\ref{corchbarra})). Note that if $(\tilde{X},\tilde{f})$ is a
time-dependent section of the vector bundle $TM\times \R\to M$ then
$\frac{\partial (\tilde{X},\tilde{f})}{\partial t}$ is the
time-dependent section given by $([\frac{\partial}{\partial
t},\tilde{X}],\frac{\partial \tilde{f}}{\partial t})$.

Now, if $\tilde{\alpha}$ is a time-dependent 1-form on $M$,
$\tilde{X}$ is a time-dependent vector field and $(x_0,t_0)\in
M\times \R$ then, using the isomorphism $T^\ast_{(x_0,t_0)}(M\times
\R)\cong T^\ast _{x_0}M\oplus T^\ast _{t_0}\R$, it follows that
\begin{equation}\label{unacuenta}
(\tilde{{\cal L}}_{\tilde{X}}\tilde{\alpha
})_{(x_0,t_0)}=({\cal L}_{\tilde{X}_{t_0}}\tilde{\alpha
}_{t_0})_{(x_0)}+\tilde{\alpha}_{(x_0,t_0)}\Big
(\frac{\partial \tilde{X}}{\partial t}_{|(x_0,t_0)}\Big ) dt_{|t_0},
\end{equation}
where $\tilde{{\cal L}}$ (resp. ${\cal L}$) is the Lie derivative on
$M\times \R$ (resp. $M$).

Moreover, $\tilde{{\cal L}}_{\frac{\partial}{\partial
t}}\tilde{\alpha }$ is a time-dependent 1-form on $M$ and if
$\tilde{f}\in C^\infty (M\times \R ,\R)$, then
$(\tilde{\alpha},\tilde{f})$ is a time-dependent section of the
vector bundle $T^\ast M\times \R\to M$ and
\begin{equation}\label{otracuenta}
\frac{\partial (\tilde{\alpha},\tilde{f})}{\partial t}=(\tilde{{\cal
L}}_{\frac{\partial}{\partial t}}\tilde{\alpha },\frac{\partial
\tilde{f}}{\partial t}).
\end{equation}
A long computation, using (\ref{ecjacobi}), (\ref{corchpoisson}),
(\ref{unacuenta}) and (\ref{otracuenta}), shows that
$$
\begin{array}{ccl}
\lcf (\tilde{\alpha},\tilde{f}),(\tilde{\beta},\tilde{g}) \rcf ^{
\, \hat{ }\, X_0}_{(\Lambda ,E)}&=&\lcf \Phi
^\ast (\tilde{\alpha}+\tilde{f}\, dt),\Phi ^\ast
(\tilde{\beta}+\tilde{g}\, dt )\rcf ^{\, \hat{ }\, X_0}_{(\Lambda ,E)}=
\lcf \tilde{\alpha}+\tilde{f}\, dt,
\tilde{\beta}+\tilde{g}\, dt \rcf _{\tilde{\Lambda}}\\
\widehat{\widetilde{\#}_{(\Lambda
,E)}}\, ^{X_0}(\tilde{\alpha},\tilde{f})&=&\widehat{\widetilde{\#}_{(\Lambda
,E)}}\, ^{X_0}(\Phi ^\ast (\tilde{\alpha}+\tilde{f}\, dt))=
\#_{\tilde{\Lambda}}(\tilde{\alpha}+\tilde{f}\, dt),
\end{array}
$$
for $\tilde{\alpha},\tilde{\beta}$ time-dependent 1-forms on $M$ and
$\tilde{f},\tilde{g}\in C^\infty (M\times \R ,\R)$, where $\tilde{\Lambda}$ is
the Poissonization of the Jacobi structure $(\Lambda ,E)$ and $\Phi ^\ast
:T^\ast (M\times \R) \to \tilde{A}^\ast =A^\ast \times \R $ is the
adjoint isomorphism of $\Phi$.

Therefore, $\tilde{A}^\ast =A^\ast \times \R$ can be
identified with $T^\ast (M\times \R)$ and, under this identification,
the Lie algebroid structure $(\lcf \, ,\, \rcf
^{\, \hat{ }\,X_0}_{(\Lambda ,E)}
,\widehat{\widetilde{\#}_{(\Lambda ,E)}}\, ^{X_0})$ is just the Lie
algebroid structure $(\lcf \, ,\, \rcf
_{\tilde{\Lambda}},\#_{\tilde{\Lambda}})$ on $T^\ast (M\times \R)$.

Consequently, using Proposition \ref{unaprop}, we deduce that, for
this particular case,  the pair
$(\tilde{A},\tilde{A}^\ast )$ is a Lie bialgebroid, when we consider
on $\tilde{A}$ and $\tilde{A}^\ast$ the Lie algebroid structures
$(\makebox{{\bf [}\, ,\, {\bf ]}}\, \bar{ }\, ^{\phi _0},
\bar{\pi}^{\phi _0})$ and $(\lcf \, ,\, \rcf
^{\, \hat{ }\,X_0}_{(\Lambda ,E)}
,\widehat{\widetilde{\#}_{(\Lambda ,E)}}\, ^{X_0})$, respectively.

In this section, we generalize the above result for an arbitrary
generalized Lie bialgebroid. In fact, we prove
\begin{theorem}\label{bialgebrizacion}
Let $((A,\phi _0),(A^\ast ,X_0))$ be a generalized Lie bialgebroid
and $(\Lambda ,E)$ the induced Jacobi structure over $M$. Consider on
$\tilde{A}$ (resp. $\tilde{A}^\ast$) the Lie algebroid structure
$(\lcf \, ,\, \rcf \, \bar{ }\, ^{\phi _0},\bar{\rho}^{\phi _0})$
(resp. $(\lcf \, ,\, \rcf _\ast
\kern-3pt \hat{ }^{X_0} ,\widehat{\rho_ \ast}^{X_0})$). Then:
\begin{itemize}
\item[{\it i)}] The pair $(\tilde{A},\tilde{A}^\ast)$ is a Lie bialgebroid over
$M\times \R$.
\item[{\it ii)}] If $\tilde{\Lambda}$ is the induced Poisson
structure on $M\times \R$ then $\tilde{\Lambda}$ is the
Poissonization of the Jacobi structure $(\Lambda ,E)$.
\end{itemize}
\end{theorem}
\prueba {\em i)} Using (\ref{primerlema}) and (\ref{corchbarra}), we
obtain that
$$
\begin{array}{ccl}
\widehat{d_\ast}^{X_0} \lcf \tilde{X},\tilde{Y}\rcf \, \bar{ }\,
^{\phi _0}&=&e^{-t}\Big ( d_\ast{}_{X_0} \lcf
\tilde{X},\tilde{Y}\rcf+d_\ast{}_{X_0}(\phi
_0(\tilde{X})\frac{\partial \tilde{Y}}{\partial t})
-d_\ast{}_{X_0}(\phi _0(\tilde{Y})\frac{\partial \tilde{X}}{\partial t})\\
&&\quad +X_0\wedge \frac{\partial}{\partial t}\lcf
\tilde{X},\tilde{Y}\rcf +X_0\wedge
\frac{\partial}{\partial t}(\phi _0(\tilde{X})\frac{\partial
\tilde{Y}}{\partial t})-X_0\wedge \frac{\partial}{\partial t}(\phi
_0(\tilde{Y})\frac{\partial \tilde{X}}{\partial t})\Big ).
\end{array}
$$
Moreover, applying (\ref{la2}), (\ref{tercerlema0}), (\ref{tercerlema1}) and
(\ref{tercerlema2}), it follows that
$$
\begin{array}{ccl}
\widehat{d_\ast}^{X_0} \lcf \tilde{X},\tilde{Y}\rcf \, \bar{ }\,
^{\phi _0}&=&e^{-t}\Big (
d_\ast{}_{X_0}\lcf \tilde{X},\tilde{Y}\rcf +X_0\wedge \lcf
\frac{\partial \tilde{X}}{\partial t},\tilde{Y}\rcf
+X_0\wedge \lcf \tilde{X},\frac{\partial \tilde{Y}}{\partial t}\rcf +\phi
_0(\tilde{X})\frac{\partial}{\partial t}(d_\ast{}_{X_0}\tilde{Y})\\
&&\quad \kern-1pt -\phi _0(\tilde{Y})\frac{\partial}{\partial
t}(d_\ast{}_{X_0}\tilde{X}) -\phi _0(\tilde{X})X_0\wedge
\frac{\partial \tilde{Y}}{\partial t}+\phi _0(\tilde{Y})X_0\kern-1pt
 \wedge \frac{\partial \tilde{X}}{\partial t}\\
&&\quad +\phi _0(\tilde{X})X_0\wedge \frac{\partial ^2\tilde{Y}} {\partial t^2}
 -\phi _0(\tilde{Y})X_0\wedge \frac{\partial
^2\tilde{X}}{\partial t^2} +\frac{\partial}{\partial t}(\phi
_0(\tilde{X}))X_0\wedge \frac{\partial \tilde{Y}}{\partial t} \\
&&\quad -\frac{\partial}{\partial t}(\phi _0(\tilde{Y}))X_0\wedge
\frac{\partial \tilde{X}}{\partial t}
+d_\ast{}_{X_0}(\phi _0(\tilde{X}))\wedge
\frac{\partial \tilde{Y}}{\partial t}
 -d_\ast{}_{X_0}(\phi _0(\tilde{Y}))\wedge \frac{\partial
\tilde{X}}{\partial t}\Big ).
\end{array}
$$
On the other hand, from (\ref{cond1}), (\ref{tercerlema0}),
(\ref{primerlema}) and (\ref{segundolema}), we have
$$
\begin{array}{ccl}
\lcf \tilde{X},\widehat{d_\ast}^{X_0} \tilde{Y}\rcf \, \bar{ }\,
^{\phi _0}&=&
e^{-t}\Big ( \lcf \tilde{X},d_\ast{}_{X_0}\tilde{Y}\rcf _{\phi _0}+
\lcf \tilde{X},X_0\wedge \frac{\partial \tilde{Y}}{\partial t}\rcf _{\phi _0}+
\phi _0(\tilde{X})\frac{\partial}{\partial t}(d_\ast{}_{X_0}\tilde{Y})\\
&&\quad +\phi _0(\tilde{X})X_0\wedge \frac{\partial
^2\tilde{Y}}{\partial t^2}-\frac{\partial \tilde{X}}{\partial t}
\wedge i_{\phi _0}(d_\ast{}_{X_0}\tilde{Y})+\phi _0(\frac{\partial
\tilde{Y}}{\partial t}) \frac{\partial \tilde{X}}{\partial t}\wedge
X_0\Big ) .
\end{array}
$$
Therefore, using (\ref{coinciden}) and (\ref{derivacion}) and the
fact that $\frac{\partial}{\partial t}(\phi _0
(\tilde{Z}))=\phi _0(\frac{\partial \tilde{Z}}{\partial t})$, for
$\tilde{Z}\in \Gamma (\tilde{A})$,
$$
\begin{array}{ccl}
\lcf \tilde{X},\widehat{d_\ast}^{X_0} \tilde{Y}\rcf \, \bar{ }^{\phi
_0} &=&e^{-t}\Big ( \lcf
\tilde{X},d_\ast{}_{X_0}\tilde{Y}\rcf _{\phi _0} +\lcf \tilde{X},X_0\rcf \wedge
\frac{\partial \tilde{Y}}{\partial t}+X_0\wedge \lcf
\tilde{X},\frac{\partial \tilde{Y}}{\partial t}\rcf
-\phi _0(\tilde{X})X_0\wedge \frac{\partial \tilde{Y}}{\partial t}\\
&&\quad +\phi_0(\tilde{X})\frac{\partial}{\partial
t}(d_\ast{}_{X_0}\tilde{Y})+\phi _0(\tilde{X})X_0\wedge
\frac{\partial ^2\tilde{Y}}{\partial t^2}-\frac{\partial
\tilde{X}}{\partial t}\wedge i_{\phi _0}(d_\ast{}_{X_0}\tilde{Y})\\
&&\quad +\frac{\partial}{\partial t}(\phi
_0(\tilde{Y}))\frac{\partial \tilde{X}}{\partial t}\wedge X_0\Big ).
\end{array}
$$
Finally, from (\ref{condcomp1}) and (\ref{cond3}), we deduce that
$$\widehat{d_\ast}^{X_0} \lcf \tilde{X},\tilde{Y}\rcf \, \bar{ }\,
^{\phi _0}=\lcf
\tilde{X},\widehat{d_\ast}^{X_0} \tilde{Y}\rcf \, \bar{ }\, ^{\phi _0}-\lcf
\tilde{Y},\widehat{d_\ast}^{X_0} \tilde{X}\rcf \, \bar{ }\, ^{\phi _0}.$$

{\em ii)} Using (\ref{cond1}), (\ref{difcomp}), (\ref{difcomp2}) and
Theorem \ref{unteorema}, we obtain that the induced Poisson structure
$\tilde{\Lambda}$ on $M\times \R$ is given by
$$
\tilde{\Lambda}(\delta\tilde{f},\delta\tilde{g})=\widehat{d_\ast}^{X_0}
\tilde{f}\cdot \bar{d}^{\phi _0}\tilde{g}=e^{-t}\Big ( d\tilde{g} \cdot d_\ast \tilde{f}+
\frac{\partial \tilde{f}}{\partial t}\rho (X_0)(\tilde{g})
-\frac{\partial \tilde{g}}{\partial t}\rho (X_0)(\tilde{f})\Big ),
$$
for $\tilde{f},\tilde{g}\in C^\infty (M\times \R ,\R)$.

On the other hand, using (\ref{Lambda}) and (\ref{E}), we prove that
$$
e^{-t}\Big ( \Lambda +\frac{\partial}{\partial t}\wedge E
\Big )(\tilde{f},\tilde{g})=
e^{-t}\Big ( d\tilde{g}\cdot d_\ast \tilde{f}+\frac{\partial
\tilde{f}}{\partial t}\rho (X_0)(\tilde{g})-\frac{\partial
\tilde{g}}{\partial t}\rho (X_0)(\tilde{f})\Big ),
$$
for $\tilde{f},\tilde{g}\in C^\infty (M\times \R ,\R)$. Therefore,
$\tilde{\Lambda}$ is the Poissonization of $(\Lambda ,E)$.\QED

Now, we discuss a converse of Theorem \ref{bialgebrizacion}.
\begin{theorem}\label{recbialg}
Let $(A,\lcf \, , \, \rcf, \rho )$ be a Lie algebroid and $\phi _0\in
\Gamma (A^\ast )$ a 1-cocycle. Suppose that $(\lcf \, , \, \rcf_\ast
, \rho _\ast )$ is a Lie algebroid structure on $A^\ast$ and that
$X_0\in \Gamma (A)$ is a 1-cocycle. Consider on $\tilde{A}=A\times
\R$ (resp. $\tilde{A}^\ast=A^\ast \times \R$) the Lie algebroid
structure $(\lcf ,\rcf \, \bar{ }\, ^{\phi _0},\bar{\rho}^{\phi _0})$
(resp. $(\lcf ,\rcf _\ast \kern-3pt \hat{ }^{X_0}
,\widehat{\rho_\ast}^{X_0})$). If $(\tilde{A},
\tilde{A}^\ast )$ is a Lie bialgebroid then the pair $((A,\phi
_0),(A^\ast ,X_0))$ is a generalized Lie bialgebroid.
\end{theorem}
\prueba Let $\{ \, ,\, \} \, \tilde{ }$ be the induced Poisson
bracket on $M\times \R$. Then, from (\ref{difcomp}), (\ref{difcomp2}) and
Theorem \ref{unteorema}, it follows that
$$\{ \tilde{f},\tilde{g} \} \, \tilde{ }=e^{-t}\Big ( d\tilde{g}\cdot d_\ast
\tilde{f}+\frac{\partial \tilde{f}}{\partial t}\rho (X_0)(\tilde{g})+
\frac{\partial \tilde{g}}{\partial t}\rho _\ast(\phi _0)(\tilde{f})
+\frac{\partial \tilde{g}}{\partial t}\frac{\partial \tilde{f}}{\partial t}
\phi _0(X_0)\Big ),$$
for $\tilde{f},\tilde{g}\in C^\infty (M\times \R,\R)$. Since $\{
\, ,\, \} \, \tilde{ }$ is skew-symmetric, we have that $\{ t,t
\} \, \tilde{ }=0$ which implies that $\phi _0(X_0)=0$. As a consequence,
$$\{ \tilde{f},\tilde{g} \} \, \tilde{ }=e^{-t}\Big (  d\tilde{g}\cdot d_\ast
\tilde{f}+\frac{\partial \tilde{f}}{\partial t}\rho (X_0)(\tilde{g})+
\frac{\partial \tilde{g}}{\partial t}\rho _\ast(\phi _0)(\tilde{f})
\Big ).$$
In particular, if $f\in C^\infty (M,\R)$ then, using that
$\{ f,t\}\, \tilde{ }=-\{ t,f\}\, \tilde{ }$, we conclude that $\rho
(X_0)=-\rho _\ast (\phi _0)$.

Now, if $X,Y\in \Gamma (A)$, from (\ref{primerlema}),
(\ref{corchbarra}) and (\ref{segundolema}), we obtain that
        $$\widehat{d_\ast}^{X_0} \lcf X,Y\rcf \, \bar{ }\, ^{\phi
               _0}=e^{-t}d_\ast{}_{X_0} \lcf X,Y\rcf ,$$
$$
\begin{array}{ccl}
\lcf X,\widehat{d_\ast}^{X_0} Y\rcf \, \bar{ }\, ^{\phi _0}-\lcf
Y,\widehat{d_\ast}^{X_0} X\rcf \, \bar{ }\, ^{\phi _0}
&=&\bar{\rho}^{\phi _0}(X)(e^{-t})d_\ast{}_{X_0}Y+e^{-t}\lcf
X,d_\ast{}_{X_0}Y\rcf \, \bar{ }\, ^{\phi _0}\\
&&-\bar{\rho}^{\phi _0}(Y)(e^{-t})d_\ast{}_{X_0}X-e^{-t}
\lcf Y,d_\ast{}_{X_0}X\rcf \, \bar{ }\, ^{\phi _0}\\
&=&e^{-t}\Big ( \lcf X,d_\ast{}_{X_0}Y\rcf _{\phi _0}-\lcf
Y,d_\ast{}_{X_0}X\rcf _{\phi _0}\Big ) .
\end{array}
$$
Thus, since
$\widehat{d_\ast}^{X_0} \lcf X,Y\rcf \, \bar{ }\, ^{\phi _0}=\lcf
X,\widehat{d_\ast}^{X_0} Y\rcf \, \bar{ }\, ^{\phi _0}-\lcf
Y,\widehat{d_\ast}^{X_0}  X\rcf \, \bar{ }\, ^{\phi _0}$, we deduce
(\ref{condcomp1}).

Finally, if $X\in \Gamma (A)$ then, using the computations in the proof of
Theorem \ref{bialgebrizacion} and the fact that
$$\widehat{d_\ast}^{X_0} \lcf X,e^tY\rcf \, \bar{ }\, ^{\phi _0}=\lcf
X,\widehat{d_\ast}^{X_0} (e^tY)\rcf \, \bar{ }\, ^{\phi _0}-\lcf
e^tY,\widehat{d_\ast}^{X_0}  X\rcf \, \bar{ }\, ^{\phi _0},$$
for all $Y\in \Gamma (A)$, we prove that $\Big ( \lcf X_0,X\rcf
+({\cal L}_\ast{}_{X_0})_{\phi _0}X\Big ) \wedge Y=0$. But this implies that
        $$\lcf X_0,X\rcf +({\cal L}_\ast{}_{X_0})_{\phi _0}X=0.$$

\vspace{-.75cm}
\QED

In \cite{MX} it was proved that if the pair $(A,A^\ast )$ is a Lie
bialgebroid then the pair $(A^\ast, A)$ is also a Lie bialgebroid.
Using this fact, Propositions \ref{isomorf} and \ref{unaprop} and
Theorems \ref{bialgebrizacion} and \ref{recbialg}, we conclude that a
similar result holds for generalized Lie bialgebroids.
\begin{theorem}
If $((A,\phi _0),(A^\ast ,X_0))$ is a generalized Lie bialgebroid, so
is $((A^\ast ,X _0),$ $(A,\phi _0))$.
\end{theorem}
\section{Triangular generalized Lie bialgebroids}
\setcounter{equation}{0}
Let $(A,\lcf \, ,\, \rcf ,\rho )$ be a Lie algebroid over $M$ and $\phi _0\in
\Gamma (A^\ast )$ a 1-cocycle. Moreover, let $P \in \Gamma
(\wedge ^2A)$ be a bisection satisfying
        $$\lcf P ,P \rcf _{\phi _0}=0.$$
We shall discuss what happens on the dual bundle $A^\ast \to M$.
Remark \ref{unasnotas} {\it i)} and Remark \ref{unasnotas2} {\it i)} suggest
us to introduce the bracket $\lcf \, ,\, \rcf _{\ast P}$ on
$\Gamma (A^\ast )$ defined by
\begin{equation}\label{corchetedual}
\begin{array}{ccl}
\lcf \phi ,\psi\rcf _{\ast P}&=&({\cal L}_{\phi _0})_{\#
_P  (\phi )}\psi - ({\cal L}_{\phi _0})_{\# _P (\psi
)}\phi-d_{\phi _0}(P (\phi ,\psi ))\\
 &=&i_{\# _P (\phi )}d_{\phi _0}\psi
-i_{\# _P (\psi )}d_{\phi _0}\phi+d_{\phi
_0}(P (\phi ,\psi )),
\end{array}
\end{equation}
for $\phi ,\psi \in \Gamma(A^\ast )$.
\begin{theorem}\label{triangulares}
Let $(A,\lcf \, ,\, \rcf ,\rho )$ a Lie algebroid over $M$, $\phi _0\in \Gamma
(A^\ast)$ a 1-cocycle and $P \in \Gamma (\wedge ^2A)$ a
bisection of $A\to M$ satisfying $\lcf P ,P \rcf _{\phi _0}=0$.
Then:
\begin{itemize}
\item[{\it i)}] The dual bundle $A^\ast \to M$ together with the bracket
defined in  (\ref{corchetedual}) and the bundle map $\rho _{\ast P} =\rho
\circ \# _P :A^\ast \to TM$ is a Lie algebroid.
\item[{\it ii)}] $X_0=-\# _P  (\phi _0)\in \Gamma(A)$ is a 1-cocycle.
\item[{\it iii)}] The pair $((A,\phi _0),(A^\ast ,X_0))$ is a
generalized Lie bialgebroid.
\end{itemize}
\end{theorem}
\prueba If we consider the Lie algebroid structure $(\lcf
,\rcf \, \bar{ }\, ^{\phi _0},\bar{\rho}^{\phi _0})$ on
$\tilde{A}=A\times \R \to M\times \R$
and the bisection $\tilde{P}=e^{-t}P \in \Gamma (\wedge
^2\tilde{A})$ then, from (\ref{corchbarra}), (\ref{difcomp2}) and
Theorem \ref{teoremaschou}, it follows that $\lcf
\tilde{P},\tilde{P}\rcf \, \bar{ }\, ^{\phi _0}=0.$
Thus, using the results of Mackenzie and Xu \cite{MX}, we deduce that
the vector bundle $\tilde{A}^\ast \to M\times \R$ admits a Lie
algebroid structure $(\lcf \, ,\, \rcf \, \tilde{ },\tilde{\rho})$, where
$\tilde{\rho}:\tilde{A}^\ast \to T(M\times \R)$ is the bundle map given by
$\tilde{\rho}=\bar{\rho}^{\phi _0}\circ \#
_{\tilde{P}}$ and $\lcf \, ,\, \rcf \, \tilde{ }$ is the
bracket on $\Gamma (\tilde{A}^\ast)$ defined by
$$\lcf \tilde{\phi},\tilde{\psi}\rcf \, \tilde{ }=i_{\# _{\tilde{P}}
(\tilde{\phi})}\bar{d}^{\phi _0}\tilde{\psi}-i_{\# _{\tilde{P}}
(\tilde{\psi})}\bar{d}^{\phi _0}\tilde{\phi}+\bar{d}^{\phi
_0}(\tilde{P}(\tilde{\phi}, \tilde{\psi})),$$
for $\tilde{\phi},\tilde{\psi}\in \Gamma(\tilde{A}^\ast)$.

Now, we will prove that
\begin{equation}\label{igualdad}
\lcf \tilde{\phi},\tilde{\psi}\rcf ^{\hat{ }X_0} _{\ast P} =
e^{-t}\Big( \lcf \tilde{\phi}, \tilde{\psi}\rcf _{\ast P}
+\tilde{\phi}(X_0)(\frac{\partial \tilde{\psi}}{\partial t}- \tilde{\psi})-
\tilde{\psi}(X_0)(\frac{\partial \tilde{\phi}}{\partial t}-\tilde{\phi})\Big )
=\lcf \tilde{\phi}, \tilde{\psi}\rcf \, \tilde{ }.
\end{equation}
From (\ref{la1}), (\ref{corchetedual}) and the definition of $X_0$, we
have that
$$
\begin{array}{ccl}
\lcf \tilde{\phi}, \tilde{\psi}\rcf ^{\hat{ }X _0} _{\ast P}
&=&e^{-t}\Big ( i_{\# _P (\tilde{\phi} )} d\tilde{\psi} -i_{\# _P
(\tilde{\psi} )} d\tilde{\phi} +d(P (\tilde{\phi} ,\tilde{\psi}))+P
(\tilde{\psi} ,\tilde{\phi})\phi _0\\
&&+(i_{\# _P (\tilde{\phi} )}\phi _0)\frac{\partial \tilde{\psi}}{\partial t}
-(i_{\# _P (\tilde{\psi} )}\phi _0)\frac{\partial
\tilde{\phi}}{\partial t}\Big ).
\end{array}
$$
Using the fact that
$$
\frac{\partial}{\partial t}(P (\tilde{\phi} ,\tilde{\psi} ))=P
(\frac{\partial \tilde{\phi}}{\partial t},\tilde{\psi} )+P (\tilde{\phi}
,\frac{\partial \tilde{\psi}}{\partial t}) = i_{\# _P (\tilde{\phi})}
\frac{\partial \tilde{\psi}}{\partial t}- i_{\# _P (\tilde{\psi})}
\frac{\partial \tilde{\phi}}{\partial t},
$$
we obtain that,
$$
\begin{array}{ccl}
\lcf \tilde{\phi}, \tilde{\psi}\rcf ^{\hat{ }X_0} _{\ast P}
&=&e^{-t}\Big ( i_{\# _P (\tilde{\phi} )}(d\tilde{\psi} +\phi _0\wedge
\frac{\partial \tilde{\psi}}{\partial t})-i_{\# _P
(\tilde{\psi} )}(d\tilde{\phi} +\tilde{\phi}
_0\wedge \frac{\partial \tilde{\phi}}{\partial t})+d(P
(\tilde{\phi} ,\tilde{\psi}
))\\
&&-(P (\tilde{\phi} ,\tilde{\psi}))\phi _0+\frac{\partial
}{\partial t}(P (\tilde{\phi}
,\tilde{\psi} ))\phi _0 \Big ).
\end{array}
$$
Thus, from (\ref{difcomp2}) and (\ref{difcomp22}), we deduce
(\ref{igualdad}).

On the other hand, using (\ref{corchtilde}) and (\ref{corchbarra}), it
follows that
$$\tilde{\rho}(\tilde{\phi})=e^{-t}\Big ( \rho _{\ast
P}(\tilde{\phi})+X_0(\tilde{\phi})\frac{\partial}{\partial
t}\Big )=\widehat{\rho _{\ast P}}^{X_0}(\tilde{\phi}),$$
for $\tilde{\phi}\in \Gamma (\tilde{A}^\ast)$. Therefore, from Proposition
\ref{equivalencia}, we prove {\it i)} and {\it ii)}.

Now, if we consider on $\tilde{A}$ (resp. $\tilde{A}^\ast$) the
Lie algebroid structure $(\lcf \, ,\, \rcf \, \bar{
}\, ^{\phi _0},\bar{\rho}^{\phi _0})$
(resp. $(\lcf \, ,\, \rcf ^{\hat{ }X_0}_{\ast P} ,$ $
\widehat{\rho _{\ast P}}^{X_0})$) then the pair $(\tilde{A},
\tilde{A}^\ast)$ is a Lie bialgebroid (see \cite{MX}). Consequently,
using Theorem \ref{recbialg}, we conclude that $((A,\phi _0),(A^\ast ,X_0))$ is
a generalized Lie bialgebroid. \QED

Let $(A,\lcf \, ,\, \rcf ,\rho )$ be a Lie algebroid and $\phi _0\in
\Gamma (A^\ast)$ a 1-cocycle. Suppose that $(\lcf \, ,\, \rcf
_\ast,\rho _\ast)$ is a Lie algebroid structure on $A^\ast$
and that $X_0\in \Gamma (A)$ is a 1-cocycle. Moreover,
assume that $((A,\phi _0),(A^\ast ,X_0))$ is a generalized Lie
bialgebroid. Then, the pair $((A,\phi _0),(A^\ast ,X_0))$
is said to be a {\em triangular generalized Lie bialgebroid} if there exists
$P \in \Gamma (\wedge ^2A)$ such that $\lcf P ,P
\rcf _{\phi _0}=0$ and
$$\lcf \, ,\, \rcf _\ast=\lcf \, ,\, \rcf _{\ast P}, \quad
\rho _\ast=\rho _{\ast _P},\quad X_0=-\# _P(\phi _0).$$
Note that a triangular generalized Lie bialgebroid $((A,\phi
_0),(A^\ast ,X_0))$ such that $\phi _0=0$ is just a {\em triangular Lie
bialgebroid} (see \cite{MX}). On the other hand, if $(M,\Lambda ,E)$
is a Jacobi manifold then, using Remarks \ref{unasnotas} and
\ref{unasnotas2}, we deduce that the pair $((TM\times
\R,(0,1)),(T^\ast M\times \R,$ $(-E,0)))$ is a triangular generalized
Lie bialgebroid.
\section{Generalized Lie bialgebras}
\setcounter{equation}{0}
In this Section, we will study generalized Lie bialgebroids over a point.
\begin{definition}
A generalized Lie bialgebra is a generalized Lie bialgebroid over a
point, that is, a pair $((\frak g ,\phi _0),(\frak g ^\ast ,X_0))$,
where $(\frak g ,[\, ,\, ])$ is a real Lie algebra of finite
dimension such that the dual space $\frak g ^\ast$ is also a Lie
algebra with Lie bracket $[\, ,\, ]_\ast$, $X_0 \in \frak g$ and
$\phi _0\in \frak g^\ast$ are 1-cocycles on $\frak g ^\ast$ and
$\frak g$, respectively, and
\begin{equation}\label{condalg1}
d_\ast {}_{X_0}[X,Y]=[X,d_\ast{}_{X_0}Y]_{\phi
_0}-[Y,d_\ast{}_{X_0}X]_{\phi _0},
\end{equation}
\begin{equation}\label{condalg2}
\phi _0 (X_0)=0,
\end{equation}
\begin{equation}\label{condalg3}
i_{\phi _0}(d_\ast{}X)+[X_0,X]=0,
\end{equation}
for all $X,Y\in \frak g$, $d_\ast$ being the differential on $(\frak
g ^\ast ,[\, ,\, ]_\ast )$.
\end{definition}
\begin{remark}
{\rm
In the particular case when $\phi _0=0$ and $X_0=0$, we recover the
concept of a {\em Lie bialgebra}, that is, a dual pair $(\frak g,\frak g
^\ast)$ of Lie algebras such that
$$ d_\ast{}[ X ,Y ] = [ X ,d_\ast Y ] -[ Y ,d_\ast X ] ,$$
for $X,Y\in \frak g$ (see \cite{D}).
}
\end{remark}
Next, we give different methods to obtain generalized Lie bialgebras.
\begin{proposition}\label{YBeq}
Let $(\frak h ,[\, ,\, ]_{\frak h})$ be a Lie algebra, $r\in \wedge ^2\frak h$
and $\bar{X}_0\in \frak h$ such that
\begin{equation}\label{YB1}
[r,r]_{\frak h}-2\bar{X}_0\wedge r=0,\qquad
[\bar{X}_0,r]_{\frak h}=0.
\end{equation}
Then, if $\frak g =\frak h \times \R$, the pair $((\frak g,
(0,1)),(\frak g ^\ast ,(-\bar{X}_0,0)))$ is a generalized Lie bialgebra.
\end{proposition}
\prueba Consider on $\frak g$ the Lie bracket $[\, ,\, ]$ given by
\begin{equation}\label{corchetedirecto}
          [(X,\lambda ),(Y,\mu )]=([X,Y]_\frak h,0),
\end{equation}
for $(X,\lambda ),(Y,\mu )\in \frak g$. One easily follows
that $\phi _0=(0,1)\in \frak h ^\ast \times \R \cong \frak g ^\ast$ is a
1-cocycle. On the other hand, the space $\wedge ^2 \frak g =\wedge
^2(\frak h \times \R )$ can be
identified with the product $\wedge ^2\frak h\times \frak h$ (see
Section \ref{algebroides}) and,
using (\ref{contraccion}), (\ref{YB1}),
(\ref{corchetedirecto}) and Theorem \ref{teoremaschou}, we have that
$P =(r,\bar{X}_0)\in \wedge ^2\frak h \times \frak h\cong \wedge
^2\frak g$ satisfies
                   $$[P ,P ]_{\phi _0}=0.$$
Therefore, from Theorem \ref{triangulares}, we deduce that there
exists a Lie bracket on $\frak g ^\ast$ and the pair $((\frak
g,(0,1)),(\frak g ^\ast ,-\# _P (0,1))$ is a
generalized Lie bialgebra. Moreover, we have that $\# _P
(0,1)=(\bar{X}_0,0)$ (see (\ref{contraccion}) and (\ref{sostenido})).\QED
\begin{remark}
{\rm
\begin{itemize}
\item[{\it i)}] Note that this method of finding generalized Lie bialgebras is
related to find algebraic Jacobi structures.
\item[{\it ii)}] Using (\ref{identificaciones}), (\ref{contraccion})
and (\ref{corchetedual}), it follows that the Lie bracket $[\, ,\,
]_\ast$ on $\frak g ^\ast$ is given by
\begin{equation}\label{YBcorch}
\begin{array}{ccl}
[(\alpha ,\lambda ),(\beta ,\mu )]_\ast &=&(coad _{\# _r (\alpha )}\beta -
coad _{\# _r (\beta )}\alpha -i_{\bar{X}_0}(\alpha \wedge \beta )\\
&&-\mu \, coad _{\bar{X}_0}\alpha +\lambda \, coad _{\bar{X}_0}\beta
,-r(\alpha ,\beta )),
\end{array}
\end{equation}
for $(\alpha ,\lambda ),(\beta ,\mu )\in \frak g ^\ast$, where $coad
:\frak h \times \frak h ^\ast \to \frak h ^\ast$ is the coadjoint
representation of $\frak h$ over $\frak h ^\ast$ defined by $(coad _X
\alpha )(Y)=-\alpha [X,Y]$, for $X,Y\in \frak h$ and $\alpha \in
\frak h ^\ast$.
\end{itemize}
}
\end{remark}
\begin{corollary}\label{YBeq2}
Let $(\frak h,[\, ,\, ]_{\frak h})$ be a Lie algebra and ${\cal
Z}(\frak h )$ the center of $\frak h$. If $r\in \wedge ^2\frak h$,
$\bar{X}_0\in {\cal Z}(\frak h)$ and
                    $$[r,r]_{\frak h}-2\bar{X}_0\wedge r=0,$$
then the pair $((\frak h ,0),(\frak h ^\ast ,-\bar{X}_0))$ is a generalized Lie
bialgebra.
\end{corollary}
\prueba Using Proposition \ref{YBeq}, we have that $((\frak g =\frak h
\times \R ,(0,1)),(\frak g ^\ast =\frak h ^\ast \times \R
,(-\bar{X}_0,0)))$ is a generalized Lie bialgebra. Furthermore, from
(\ref{YBcorch}) and since $\bar{X}_0\in {\cal Z}(\frak h)$, we deduce
that the Lie bracket on $\frak
g ^\ast$ is given by
\begin{equation}\label{YBcorch2}
[(\alpha ,\lambda ),(\beta ,\mu )]_\ast =(coad _{\# _r (\alpha )}\beta -
coad _{\# _r (\beta )}\alpha -i_{\bar{X}_0}(\alpha \wedge \beta ),-r(\alpha
,\beta )),
\end{equation}
for $(\alpha ,\lambda ),(\beta ,\mu )\in \frak g ^\ast$. Then $\frak
h$ and $\frak h ^\ast$ are Lie algebras, where the Lie bracket on
$\frak h ^\ast$ is defined by
\begin{equation}\label{YBcorch3}
[\alpha ,\beta ]_{\frak h ^\ast}=coad _{\# _r (\alpha )}\beta -
coad _{\# _r (\beta )}\alpha -i_{\bar{X}_0}(\alpha \wedge \beta ),
\end{equation}
for $\alpha ,\beta \in \frak h ^\ast$. Moreover, using
(\ref{corchetedirecto}), (\ref{YBcorch2}), (\ref{YBcorch3}) and
the fact that $((\frak g ,(0,1)),(\frak g ^\ast ,(-\bar{X}_0,0)))$ is a
generalized Lie bialgebra, we conclude that $((\frak h ,0),(\frak h
^\ast ,-\bar{X}_0))$ is a generalized Lie bialgebra.\QED
\begin{remark}
{\rm From (\ref{YBcorch3}) and Corollary \ref{YBeq2}, we deduce a
well-known result (see \cite{D}): if $(\frak h ,[\, ,\, ]_{\frak
h})$ is a Lie algebra, $r\in \wedge ^2 \frak h$ is a solution of the
classical Yang-Baxter equation (that is, $[r,r]_{\frak h}=0$) and on
$\frak h ^\ast$ we consider the bracket defined by
$$[\alpha ,\beta ]_{\frak h ^\ast}=coad _{\# _r (\alpha )}\beta -
coad _{\# _r (\beta )}\alpha , \mbox{  for }\alpha ,\beta \in \frak h ^\ast ,$$
then the pair $(\frak h ,\frak h ^\ast )$ is a Lie bialgebra.
}
\end{remark}
\begin{examples}
\begin{itemize}
{\rm
\item[{\it i)}] Let $(\frak h ,[\, ,\, ]_{\frak h})$ be the Lie algebra of the
Heisenberg group $H(1,1)$. Then, $\frak h =<\{ e_1,e_2,e_3\}>$ and
$$[e_1,e_2]_{\frak h} =e_3, \quad e_3\in {\cal Z}(\frak h).$$
If we take $r=e_1\wedge e_2$ and $\bar{X}_0=-e_3$, we have that
$[r,r]_{\frak h}-2\bar{X}_0\wedge r=0$ and thus, using Corollary
\ref{YBeq2}, we conclude that $((\frak h,0),(\frak h ^\ast ,e_3))$ is a
generalized Lie bialgebra. Note that $r$ and $\bar{X}_0$ induce the
canonical left-invariant contact structure of $H(1,1)$.
\item[{\it ii)}] Denote by $\frak h =\frak s \frak u (2)=\{ A\in gl (2, \C )/
\bar{A}^T=-A, trace A=0\}$ the Lie algebra of the special unitary
group $SU(2)$ and by $\sigma
_1,\sigma _2$ and $\sigma _3$ the Pauli matrices
$$\sigma _1 =\left ( \begin{array}{cc} 0&1\\1&0\end{array} \right )
\quad
\sigma _2 =\left ( \begin{array}{cc} 0&-i\\i&0\end{array} \right )
\quad
\sigma _3 =\left ( \begin{array}{cc} 1&0\\0&-1\end{array} \right ).
$$
Then, the matrices $\{ e_1=\frac{i}{2}\sigma _1,e_2=\frac{i}{2}\sigma
_2, e_3=\frac{i}{2}\sigma _3\}$ form a basis of $\frak h =\frak s
\frak u (2)$ and if $[\, ,\, ]_{\frak h}$ is the Lie bracket on
$\frak h$, we have that
$$[e_1,e_2]_{\frak h}=-e_3, \quad [e_1,e_3]_{\frak h}=e_2,\quad
[e_2,e_3]_{\frak h}=-e_1.$$
Since $r=e_1\wedge e_2$ and $\bar{X}_0=e_3$ satisfy the equations
(\ref{YB1}), $((\frak s \frak u (2)
\times \R ,(0,1)),(\frak s \frak u (2) ^\ast \times \R ,(-\bar{X}_0,0)))$ is a
generalized Lie bialgebra. Note that if $[\, ,\, ]$ is the Lie
bracket on $\frak s \frak u (2)\times \R$ then, from
(\ref{corchetedirecto}), it follows that  $(\frak s \frak u (2)\times
\R,[\, ,\, ])$ is just the Lie algebra of the unitary group $U(2)$.
Moreover, $r$ and $\bar{X}_0$ induce a left-invariant contact
structure on $SU(2)$.
\item[{\it iii)}] Let $Gl(2,\R)$ be the general linear group and
$\frak h = gl(2,\R )$ its Lie algebra. A basis of $\frak h$ is given
by the following matrices
$$e_1 =\left ( \begin{array}{cc} 0&1\\0&0\end{array} \right )
\quad
e_2 =\left ( \begin{array}{cc} 0&0\\1&0\end{array} \right )
\quad
e_3 =\left ( \begin{array}{cc} 1&0\\0&-1\end{array} \right )
\quad
e_4 =\left ( \begin{array}{cc} 1&0\\0&1\end{array} \right ) .
$$
If $[\, ,\, ]_{\frak h}$ is the Lie bracket on $\frak h$, we have that
    $$[e_1,e_2]_{\frak h}=e_3,\quad [e_1,e_3]_{\frak h}=-2e_1,\quad
                       [e_2,e_3]_{\frak h}=2e_2,\quad
                  e_4\in {\cal Z}(\frak h ).$$
Therefore, if $r=e_1\wedge e_3 + (e_1-\frac{1}{2}e_3)\wedge e_4$ and
$\bar{X}_0=-e_4$, we deduce that
                   $$[r,r]_{\frak h}-2\bar{X}_0\wedge r=0.$$
Consequently,  $((\frak h ,0),(\frak h ^\ast ,-\bar{X}_0))$ is a generalized
Lie bialgebra (see Corollary \ref{YBeq2}).
}
\end{itemize}
\end{examples}

\vspace{.2cm}
{\small {\bf Acknowledgments.} Research partially supported by DGICYT
grant PB97-1487 (Spain). D. Iglesias wishes to thank Spanish Ministerio
de Educaci\'on y Cultura for a FPU grant.}

\end{document}